\documentclass[11pt]{article}

\makeatletter
\newfont{\footsc}{cmcsc10 at 8truept}
\newfont{\footbf}{cmbx10 at 8truept}
\newfont{\footrm}{cmr10 at 10truept}
\makeatother
\newcommand{\half}{{\textstyle\frac12}}

\newcommand{\threehalves}{{\textstyle\frac32}}
\newcommand{\threequarters}{{\textstyle\frac34}}
\newcommand{\threeoverseven}{{\textstyle\frac37}}
\newcommand{\quarter}{{\textstyle\frac14}}

\newcommand{\ignore}[1]{}

\newtheorem{theorem}{Theorem}
\newtheorem{lemma}{Lemma}
\newtheorem{example}{Example}

\newtheorem{corollary}{Corollary}

\newtheorem{proposition}{Proposition}

\newenvironment{proof}{\begin{trivlist}
                       \item[]{\bf Proof}
                       \hspace{0cm} }{\hfill {\large $\bullet$}
                       \end{trivlist}}

\def\quad{ ~ }
\def\qquad{ ~~ }

\input{epsf}

\def\ssp{\def\baselinestretch{1.5}\large\normalsize}
\def\dsp{\def\baselinestretch{2.0}\large\normalsize}

\title{Almost Product Evaluation of Hankel Determinants}

\author{ \ \\
{\small \"Omer E\u{g}ecio\u{g}lu }\\
{\small Department of Computer Science, }
{\small University of California,}\\
{\small Santa Barbara CA 93106 \   ({\tt omer@cs.ucsb.edu})}\\ \ \\
{\small Timothy Redmond}\\
{\small Stanford Medical Informatics, Stanford University}\\
{\small Stanford, CA 94305 \ ({\tt tredmond@stanford.edu})}\\ \ \\
{\small Charles Ryavec}\\
{\small College of Creative Studies, University of California,} \\
{\small Santa Barbara CA 93106 \   ({\tt ryavec@math.ucsb.edu})} }

\begin{document}

\maketitle

\begin{abstract}
An extensive literature exists describing various techniques for the evaluation of Hankel
determinants. The prevailing methods such as Dodgson condensation, continued fraction expansion, LU decomposition, 
all produce product formulas when they are applicable. We mention the classic case of the Hankel determinants
with binomial entries ${3 k +2 \choose k}$ 
and those with entries ${3 k \choose k}$; both of these classes of Hankel determinants 
have product form evaluations. 
The intermediate case,  ${3 k +1 \choose k}$ has not been evaluated. 
There is a good reason for this: these latter determinants do not have product form evaluations.

In this paper we evaluate the Hankel determinant of ${3 k +1 \choose k}$. The evaluation is a sum 
of a small number of products, an {\em almost product}. 
The method actually provides more, and as applications, we present the salient points for the evaluation of 
a number of other Hankel determinants with polynomial 
entries, along with product and almost product form evaluations at special points.

\end{abstract}

\ssp

\section{Introduction}
\label{intro}

A determinant 
$$ 
H_n = \det [a_{i,j}]_{0\le i,j\le n} 
$$ 
whose entries satisfy
$$ a_{i,j} = a_{i+j}$$ 
for some sequence
$ \{a_k\}_{k \geq 0}$ is said to be a {\em Hankel determinant}.
Thus $H_n$ is the determinant of a special type of $ (n+1) \times (n+1)$ symmetric matrix.

In various cases of  Hankel determinant evaluations, special techniques 
such as Dodgson condensation,  continued fraction expansion, and LU decomposition are applicable. These methods 
provide product formulas for a large class of Hankel determinants.
A modern treatment of the theory of determinant evaluation including Hankel determinants as well as a
substantial bibliography can be found 
in Krattenthaler \cite{K99,K05}. 

The product form determinants of special note
are those whose factors have some particular attraction.
Factorials and other familiar combinatorial entities that appear as factors have an especially pleasing quality, and 
we find an extensive literature devoted to the evaluation of classes of Hankel determinants as such products.

Several classical Hankel determinants involve entries that are binomial coefficients or
expressions closely related to
binomial coefficients.
Perhaps the most well-known of these is where $ a_k = { 2k+1 \choose k}$, 
and $ a_k = \frac{1}{2k+1} { 2k+1 \choose k}$, 
for which $ H_n =1$ for all $n$. 

It is also known that  
binomial entries such as 
\begin{equation}
\label{exaks}
 a_k =  { 3k \choose k}, ~~
 a_k =   { 3k+2 \choose k}
\end{equation}
yield product evaluations for the corresponding $H_n$.
In fact, product formulas have been shown to 
exist for a host of other cases (see Gessel and Xin \cite{Gessel05}),
and we only mention

$$
a_k = \frac{1}{3k+1}{ 3k+1 \choose k} , ~ ~
a_k = \frac{9k+14}{(3k+4)(3k+5)}{ 3k+2 \choose k+1} 
$$
as representatives. 

However,  within the restricted class of 
Hankel determinants defined by the binomial coefficients
\begin{equation}
\label{form}
a_k = a_k^{(\beta , \alpha )} = {  \beta k + \alpha \choose k},
\end{equation}
parametrized by a pair of integers $\beta>0$ and $\alpha$, 
it is a rare phenomenon that the determinant evaluations are in product form.
An extensive check of Hankel determinants of sequences 
$a_k$ in the form  (\ref{form}) 
suggests that 
there is no product formula for $H_n$ in general for such binomial sequences.
In fact, it would seem that the instances for which $H_n$ has a product form can be enumerated in full:
\begin{enumerate}
\item[(i)] $ \beta  = 1$, $ \alpha $ arbitrary,
\item[(ii)] $  \beta = 2$, $ \alpha =0, 1, 2,3,4$,
\item[(iii)] $  \beta = 3$, $ \alpha  =0, 2$.
\end{enumerate}

All the other cases are likely not products, but in any event, the
question remains open: are evaluations 
possible in these cases?

Let $H_n^{( \beta, \alpha)}$ denote the $(n+1) \times (n+1)$ Hankel determinant with entries
$ a_k^{(\beta , \alpha)} $ as defined in (\ref{form}). We 
will also  use the term {\em $(\beta , \alpha )$-case}  to 
refer to the evaluation of 
$H_n^{( \beta, \alpha)}$.

Numerical data indicates the intriguing 
possibility that the $H_n^{( \beta, \alpha)}$ might be evaluated as a sum of a small number of
products, where ``small'' would mean $ O (n^d)$ summands for some fixed $ d = d( \beta , \alpha )$.
We refer to such an 
evaluation as an {\em almost product}.

The evidence of an almost product evaluation
of $H_n^{( \beta, \alpha)}$ is most pronounced for $\beta =  3$, 
and we begin with some sample data.

For the $ (3, 2 )$-case  the Hankel determinants evaluate to

{\footnotesize
\begin{eqnarray*}
H_{10}^{(3,2)} &= & 2^2 \cdot 3 \cdot 7^3 \cdot 37 \cdot 41^2 \cdot 43^3 \cdot 47^3 \cdot 53^2 \cdot 59  \cdot 61\\
H_{20}^{(3,2)} &=& 3^7 \cdot 11 \cdot 17 \cdot 29^2 \cdot 31 \cdot 67 \cdot 71^2 \cdot 73^3 \cdot 79^5 \cdot 83^6
\cdot 89^6 \cdot 97^5 \cdot 101^4 \cdot 103^4 \cdot 107^3 \cdot 109^3 \cdot 113^2 \\
H_{30}^{(3,2)} & = & 2^{10} \cdot 5^{12} \cdot 11^9 \cdot 13^3 \cdot 41^3 \cdot 43^3 \cdot 97 \cdot 101^2 \cdot 103^3 \cdot
107^4 \cdot 109^5 \cdot 113^6 \cdot 127^10 \cdot 131^9 \cdot 137^8 \cdot 139^8 \cdot 149^6 \cdot \\ 
&~&   151^6 \cdot 157^5 \cdot 163^4 \cdot 167^3 \cdot 173^2 \cdot 179 \cdot 181
\end{eqnarray*}
}
The small prime factors are indicative of the fact that there is an underlying product formula. In fact for
the $(3,2) $-case,  the Hankel
determinant is explicitly given by
\begin{equation}
\label{3np2}
   H_n^{(3,2)} = \prod_{i=1}^n\frac{(6i+4)! (2i+1)!}{2 (4 i+2)!(4i+3)!}  
\end{equation}

We mention also the $ (3,0) $-case
for which the Hankel determinant also has small prime factors, and possesses the  product  evaluation
\begin{equation}
\label{3np0}
   H_n^{(3,0)} = \prod_{i=1}^n\frac{3(3i+1) (6i)!(2i)!}{ (4i)!(4i+1)!}
\end{equation}

For the 
$(3,1) $-case, 
we get the following intriguing evaluations:
{\footnotesize
\begin{eqnarray*}
H_{10}^{(3,1)} &=& 2^2 \cdot 7^2 \cdot {37} \cdot {41^2} \cdot {43^3} \cdot {47^2} \cdot 53  \cdot {41740796329}\\
H_{20}^{(3,1)}  &=  &{3^ 8} \cdot {29} \cdot {67} \cdot {71^2} \cdot {73^3} \cdot {79^5} \cdot {83^6} \cdot {89^5}
\cdot {97^4} \cdot {101^3} \cdot {103^3}\cdot {107^2} \cdot {109^2} \cdot {113} \cdot {631} \cdot
{548377971864917477341}\\
H_{30}^{(3,1)} & =& {2^{10}} \cdot {5^{10}} \cdot {11^9} \cdot {13^2} \cdot {41^3} \cdot {43^2} \cdot {97} \cdot
{101^2}\cdot {103^3} \cdot {107^4} \cdot {109^5} \cdot {113^6} \cdot {127^9} \cdot {131^8} \cdot {137^7} \cdot
{139^7} \cdot {149^5} \cdot \\
&~&  {151^5} \cdot {157^4} \cdot {163^3} \cdot {167^2} \cdot {173} \cdot {569} \cdot {920397320923}
\cdot {56029201596264233691799}
\end{eqnarray*}
}
The existence of large primes in the factorizations indicates that $H_n^{(3,1)} $ 
does not have a product form evaluation. However 
if we write 
$H_n = P_n Q_n$ where $P_n$ is the product 
of the small primes and $Q_n$ is the product of the 
large primes left over, then 
it appears that the estimates $ ~ \log P_n = \Omega ( n^2)$ and $ \log Q_n = O (n) $ hold.
This suggests that these Hankel determinants can be represented 
as a sum of  $O (n)$  number of products, all of which have very similar
representations.

The purpose of this paper is to provide a method that evaluates $ H_n^{(3,1)}$ and a number of other Hankel determinants 
as almost products.  For $H_n= H_n^{(3,1)}$, we obtain

\begin{equation}
\label{mainformula1}
H_n  =
      (-1)^n \prod_{i=1}^n \frac{ (6i-3)!(3i+2)!(2i-1)!} {(4i-1)!(4i+1)!(3i-2)!}
   \sum_{i=0}^n\frac{n!(3n+i+2)! (-6)^i}
                    {(3n+2)! (n-i)! (2i+1)!},
\end{equation}
or alternately
\begin{equation}
\label{mainformula2}
H_n  =
 \prod_{i=1}^n\frac{(6i+4)! (2i+1)!}{2 (4 i+2)!(4i+3)!}
       \sum_{i=0}^n \frac{ n!(4n+3)!!  (3n+i+2)! }
                              {(3n+2)! i!(n-i)!(4n+2i+3)!!},
\end{equation}
each as a sum of $n+1$ products.
The method actually evaluates more, and we describe now the general situation in
the $(3,1)$-case, which consists of three
basic ingredients:

\begin{enumerate}
\item[(I)] Replace $a_k$ with polynomials 
\begin{equation}
\label{akintro31}
a_k(x)= a_k^{(3 , 1 )} (x) =  \sum_{m=0}^{k}{3 k+ 1 -m \choose k-m}x^{m}
\end{equation}
so that $ a_k(x)$ is a monic polynomial 
of degree $k$ with $ a_k = a_k(0)$. 
\item[(II)] Show that the $a_k(x)$ satisfy certain differential-convolution equations.
\item[(III)] Show that the resulting determinants $H_n(x)$ themselves satisfy certain differential equations. \\
\end{enumerate}

The $(n+1) \times (n+1)$ Hankel determinant 
$H_n (x) = H_n^{(3,1)}(x)$ 
is then expressed as the 
power series solution of the differential equation in (III), and we give it here as
it is stated as Theorem \ref{final3np1}:

\begin{equation}
\label{mainformula}
H_n (x)  =
      (-1)^n \prod_{i=1}^n \frac{ (6i-3)!(3i+2)!(2i-1)!} {(4i-1)!(4i+1)!(3i-2)!}
   \sum_{i=0}^n\frac{n!(3n+i+2)! 2^i(x-3)^i}
                    {(3n+2)! (n-i)! (2i+1)!}.
\end{equation}
An alternate expression for this evaluation appears in Theorem \ref{thme5} in Section \ref{additional}.

Similarly, steps (I), (II), (III), {\em mutatis mutandis}, yield
(Theorem \ref{final2np1}):
$$
    H_n^{(2,1)}(x) = (-1)^n (2n+1)
\sum_{i=0}^{n}\frac{(n+i)!2^i (x-2)^i}{(n-i)! (2i+1)!} ,
$$
where
\begin{equation}
\label{aks2np1}
a_k(x) = a_k^{(2,1)} (x)=  \sum_{m=0}^{k}{2k+1-m \choose k-m}x^{m} ~.
\end{equation}

These polynomial families have a number of interesting properties that we briefly discuss. For example, 
the polynomials $H_n^{(3,1)}(x)$ satisfy a three-term recursion, their roots are real, and interlace. Furthermore the
specializations at $ x = 3, \threehalves, \threequarters$ all have product evaluations.
The polynomials
$H_n^{(2,1)}(x)$  form an  orthogonal family. A few other classes of Hankel determinants can be evaluated in almost product
form by simple transformations of these polynomials (e.g., Example \ref{example3k0m} 
in Section \ref{additional}, and Corollaries \ref{2np2jacobi1} and \ref{2np2jacobi2} in Section
\ref{21desection}).

Returning briefly to the general case of the determinants  
 $H_n^{( \beta, \alpha)}$, we remark that there are a few more cases which fall under the method 
described above. These would also include the same three ingredients:

\begin{enumerate}
\item[(I)] Replace $a_k$ with polynomials
\begin{equation}
\label{akintro}
a_k(x)= a_k^{(\beta , \alpha )} (x) =  \sum_{m=0}^{k}{\beta  k+ \alpha  -m \choose k-m}x^{m}
\end{equation}
so that $ a_k(x)$ is a monic polynomial
of degree $k$ with $ a_k = a_k(0)$.
\item[(II)] Show that the $a_k(x)$ satisfy certain differential-convolution equations.
\item[(III)] Show that the resulting determinants $H_n(x)$ themselves satisfy certain differential equations. \\
\end{enumerate}

The $(3,1)$-case and the $(2,1)$-case are governed by second order differential equations, but even these 
cases already present considerable technical problems to overcome. 
We mention some further difficulties that arise in the consideration of other $(\beta, \alpha)$-cases 
in Section \ref{conclusions}. 
A number of additional almost product evaluations of Hankel 
determinants are given in Section \ref{additional}.
For these additional results given as Theorems \ref{thme1}, \ref{thme2}, \ref{thme3}, \ref{thme4}, \ref{thme5}
and special product evaluations that appear in (\ref{special1a}), (\ref{special1b}) and (\ref{special2}),
we provide the necessary identities for proving the differential equations, 
mimicking the proofs we present for $H_n^{(3,1)}(x)$ and  $H_n^{(2,1)}(x)$.

Finally, a remarkable property of these $(n+1) \times (n+1)$ 
Hankel determinants $H_n^{( \beta, \alpha)}(x)$ is that the degree of the polynomial $H_n(x)$ is only $n$, 
indicating an extraordinary amount of cancellation in the expansion of the determinant.
The unusual degree of cancellation is a basic property 
of a large class of Hankel determinants with polynomial entries. This class contains the 
Hankel determinants defined by 
polynomials
in (\ref{akintro}) that we consider. The degree result
is of independent interest, and we include an exact statement and a proof of it in Appendix III. \\

We would like to remark that 
the differential-convolution equations (II) used in this paper are reminiscent of 
the equations that arise in the study of the Painlev\'e II  equation and
the Toda lattice \cite{Kajiwara01,Iwasaki02}.  

\section{The $(3,1)$-case}

\subsection{Differential-convolution equations}

In the proof of the $(3,1)$-case, we denote 
the polynomials 
$a_k^{(3,1)}(x) $ by $ a_k$, $H_n^{( 3, 1)} (x) $  by $H_n$, and the differentiation operator by $d_x$.

We need two identities given below in Lemmas
\ref{lemma1} and \ref{lemma2}. 
The first is a 
differential-convolution equation.
The second identity involves convolutions and $a_k$
but no derivatives.
The 
proofs are given in Appendix II.
\begin{lemma}
\label{lemma1}
Let the polynomials $a_k = a_k^{(3,1)} (x)$ be as defined in (\ref{akintro31}). 
Then
\begin{eqnarray}
\label{FirstId31}
    (x-3)(2x-3)(4x-3)d_xa_k & = & 2(2k+3)a_{k+1} 
       -(8x^2-18x+27k+36)a_k\\ \nonumber
        & & + ~4(2x^2-6x+3)c_k    
       -27(2x^2-6x+3)c_{k-1}     
\end{eqnarray}
where
\begin{equation}
\label{conv}
c_k=   c_k (x) = \sum_{m=0}^{k}a_m(x) a_{k-m}(x) ,~~~~~~~~(c_{-1}=0).
\end{equation}
\end{lemma}
\begin{lemma}
\label{lemma2}
With $a_k= a_k^{(3,1)}(x)$ and $c_k= c_k(x)$ as in Lemma \ref{lemma1},
we have 
\begin{eqnarray*}\label{SecondId31}
    & & 4(2k+5)(x-1) a_{k+2}-(2(16x^3-72x^2+135x-81)k
                                +2(24x^3-92x^2+180x-117))a_{k+1} \\
        & & ~~~~~~~~ + ~ (27(2x-3)^3 k+54(2x-3)(2x^2-4x+3))a_k 
        + 8(x-1)(2x^2-6x+3)c_{k+1} \\
        & & ~~~~~~~~ + ~ 2(8x^4-114x^3+324x^2-297x+81)c_k 
         -27x(2x-3)(2x^2-12x+9)c_{k-1} = 0 
\end{eqnarray*}
\end{lemma}

Lemmas \ref{lemma1} and \ref{lemma2} will be needed for the proof of the differential equation
satisfied by  the determinants $H_n$. This differential equation is 
given below in Theorem \ref{thm1}.

In addition to the first two 
identities in Lemma \ref{lemma1} and \ref{lemma2}, a much more complicated third identity
involving the $a_k$ is also needed for the proof of 
this differential equation. This will emerge in the course of the proof of (\ref{de1}).

\begin{theorem}
\label{thm1}
Let the polynomials $a_k = a_k^{(3,1)} (x)$ be as in (\ref{akintro31})
and define the $(n+1) \times (n+1)$ Hankel matrix by 
$$
    A_n = A_n(x) = [a_{i+j}(x)]_{0\le i,j\le n} 
$$
Then
\begin{equation}
\label{defHn}
    H_n= H_n^{(3,1)} (x) = \det A_n (x)
\end{equation}
satisfies the differential equation
\begin{equation}
\label{de1}
      (x-1)(x-3)d_x^2y + \left( 2(n+2)(x-3)+3 \right) d_xy - 3n(n+1)y = 0.
\end{equation}
\end{theorem}
\begin{proof}
This is easy to check for $n=0, 1$. Henceforth
we assume that $n\ge 2$. 

We will find expression for the first derivative $d_x H_n$ Section 
\ref{firstdersection} and the second derivative 
$d_x^2 H_n$ Section \ref{seconddersection}. This is followed by Section \ref{evaluatex} 
on specializations of $x$: product formulas for 
$H_n(3)$, $H_n ( \threehalves )$ and $H_n( \threequarters )$ 
are given as three corollaries in 
Sections \ref{3section}, \ref{3over2section}, and \ref{3over4section}.
The specializations make use of a Dodgson-like expansion result we prove as 
Proposition \ref{detidentprop}
at the start of Section \ref{evaluatex}.
The proof of 
Theorem \ref{thm1} continues in Section \ref{lastpiece} where 
we put together the expressions obtained for the derivatives and the third identity mentioned to 
prove (\ref{de1}).

\subsection{Calculating the First Derivative}
\label{firstdersection}

The first step is to find a simple form for the derivative of $H_n$. We begin with the expression
\begin{equation}
\label{der}
     d_x H_n = \mbox{Tr}(A_n^{-1}d_xA_n)H_n
\end{equation}
for the derivative of a determinant, 
where 
$$
     d_xA_n =
     d_xA_n (x) =
    [d_x a_{i+j}(x)]_{0\le i,j\le n} .
$$

Referring to Lemma \ref{lemma1}, 
we write 
$$
    (x-3)(2x-3)(4x-3)\mbox{Tr}(A_n^{-1}d_xA_n)
$$
as 
\begin{eqnarray}
\label{term1}
      && \mbox{\ \ \ } 2\mbox{Tr}(A_n^{-1}[(2(i+j)+3)a_{i+j+1}]_{0 \leq i,j \leq n})   \\
\label{term2}
                && + \mbox{Tr}(A_n^{-1}[- (8x^2-18x+27(i+j)+36)a_{i+j}]_{0 \leq i,j\leq n})  \\
\label{term3}
                &&+ 4(2x^2-6x+3)\mbox{Tr}(A_n^{-1}[c_{i+j}]_{0 \leq i,j \leq n})     \\
\label{term4}
                &&-27(2x^2-6x+3)\mbox{Tr}(A_n^{-1}[c_{i+j-1}]_{0 \leq i,j \leq n}) 
\end{eqnarray}
where the convolutions $c_k$ are defined in (\ref{conv}).

We render each of these four expressions (\ref{term1})-(\ref{term4}) in a simple form, and then
combine them all into an expression for the derivative in (\ref{der}). After that is done, 
we go through a similar computation for the second derivative of $H_n$, where we use the recursion
in Lemma \ref{lemma2} for the simplifications. The differential equation
will follow from a third identity, the proof of which makes up the bulk of the work for the rest of the 
argument.

We begin with the trace term (\ref{term1}):
Let I denote the identity matrix of relevant dimension and define
the two matrices
\begin{eqnarray*}
    B_n = B_n(x) &  = &\left[\begin{array}{ccccc}
          a_1 & a_2 &  \ldots & a_{n+1} \\
          a_2 & a_3 &  \ldots & a_{n+2} \\
       \vdots &     &        & \vdots \\
          a_{n+1} & a_{n+2}    &  \ldots       & a_{2n+1}
       \end{array}\right] \\
     L_n & = & \left[\begin{array}{ccccc}
               0 \\
                 & 1 \\
                 &   & 2 \\
                 &   &   & \ddots \\
                 &   &   &        & n
             \end{array}\right].
\end{eqnarray*}

Then
$$
  2\mbox{Tr}\left(
               A_n^{-1} \left[(2(i+j) + 3)a_{i+j+1}(x)\right]_{0 \leq i,j \le n}
              \right) =
  2\mbox{Tr}\left(
               A_n^{-1} ((2L_n+\frac{3}{2}I)B_n + B_n(2L_n+\frac{3}{2}I))
              \right)
$$

Now define $\sigma_0, \sigma_1, \ldots,\sigma_n$ and $K_n $ as follows:
\begin{equation}
\label{sigman}
   \left[\begin{array}{c}
       \sigma_0 \\
       \sigma_1 \\
       \vdots \\
       \sigma_n
   \end{array}\right] = A_n^{-1}\left[\begin{array}{c}
                                a_{n+1} \\
                                a_{n+2} \\
                                \vdots \\
                                a_{2n+1}
                               \end{array}\right]
\end{equation}
and
\begin{equation}
\label{defKn}
   K_n = \det \left[\begin{array}{cccccc}
       a_0     & a_1     & \ldots & a_{n-1}  & a_{n+1} \\
       a_1     & a_2     & \ldots & a_{n}    & a_{n+2} \\
       \vdots  &         & \ddots                        \\
       a_{n-1} & a_n     & \ldots & a_{2n-2} & a_{2n} \\
       a_{n} & a_{n+1} & \ldots & a_{2n-1}   & a_{2n+1}
              \end{array}\right] .
\end{equation}
By Cramer's rule we have
\begin{equation}
\label{Cramersrule1}
  \sigma_n = \frac{K_n}{H_n}.
\end{equation}
Therefore 
\begin{equation}\label{A1AsSigma}
A_n^{-1}B_n = \left[\begin{array}{cccccccc}
                     0 &   &   & \ldots &   & \sigma_0 \\
                     1 & 0 &   &        &   & \sigma_1 \\
                       & 1 & 0 &        &   & \sigma_2 \\
                       &   &   & \ddots &   &  \vdots  \\
                       &   &   &        & 0 & \sigma_{n-1} \\
                       &   &   &        & 1 & \sigma_n
               \end{array}\right].
\end{equation}
and 
$$
   A_n^{-1}B_n(2L_n+\frac{3}{2}I) = \left[\begin{array}{cccccccc}
                     0   &     &   & \ldots &        & (2n+3/2)\sigma_0 \\
                     3/2 & 0   &   &        &        & (2n+3/2)\sigma_1 \\
                         & 7/2 & 0 &        &        & (2n+3/2)\sigma_2 \\
                         &     &   & \ddots &        &  \vdots  \\
                         &     &   &        & 0      & (2n+3/2)\sigma_{n-1} \\
                         &     &   &        & (4n-1)/2 & (2n+3/2)\sigma_n
               \end{array}\right].
$$
Since
$$
      (2L_n+\frac{3}{2}I)B_n A_n^{-1}  =  (A_n^{-1}B_n(2L_n+\frac{3}{2}I))^T
$$
we can write
\begin{eqnarray*}
  \lefteqn{2\mbox{Tr}\left(
               A_n^{-1} ((2L_n+\frac{3}{2}I)B_n + B_n(2L_n+\frac{3}{2}I))
              \right)} \\
   & = & 2\mbox{Tr}\left(A_n^{-1}(2L_n+\frac{3}{2}I)B_n \right) +
              2\mbox{Tr}\left(A_n^{-1} B_n (2L_n+\frac{3}{2}I)\right) \\
   & = & 2\mbox{Tr}\left((2L_n+\frac{3}{2}I) B_n A_n^{-1}\right) +
              2\mbox{Tr}\left(A_n^{-1} B_n (2L_n+\frac{3}{2}I)\right) \\
   & = & 4\mbox{Tr}\left(A_n^{-1} B_n (2L_n+\frac{3}{2}I)\right) \\
   & = & 4(2n+\frac{3}{2})\sigma_n \\
   & = & 2(4n+3)\sigma_n
\end{eqnarray*}
The last expression gives
\begin{equation}
\label{term1simplified}
       2\mbox{Tr}(A_n^{-1}[(2(i+j)+3)a_{i+j+1}]_{0 \leq i,j \leq n})   =
   2(4n+3)\frac{K_n}{H_n}.
\end{equation}
for the desired form of the first term (\ref{term1}).
It is useful to record in passing that
\begin{equation}
\label{DeltaOneAsTrace31}
    \frac{K_n}{H_n} = \sigma_n = \mbox{Tr}(A_n^{-1}B_n).
\end{equation}
The identity (\ref{DeltaOneAsTrace31}) 
will be useful later when we calculate the
derivative of $K_n$.

Now we consider the second trace term (\ref{term2}).
The calculation of this term is done in the same manner as the first
term but the evaluation is somewhat simpler.  We get the expression:
\begin{equation}
\label{term2simplified}
  \mbox{Tr}(A_n^{-1}[ - (8x^2-18x+27(i+j)+35)a_{i+j}]_{0 \leq i,j \leq n})=
   -(n+1)(8x^2-18x+27n+36).
\end{equation}

The final two terms 
(\ref{term3}) and (\ref{term4}) 
require a new technique. 
We will use ideas from \cite{Kajiwara01,Iwasaki02}
where an identity similar to the following is used:
\begin{equation}
\label{M0Equation}
  \left[c_{i+j}\right]_{0 \leq i,j \leq n} = E_n A_n + A_n E_n^T
\end{equation}
where
\begin{equation}
\label{Ematrix}
   E_n= E_n(x) = \left[\begin{array}{ccccc}
      a_0/2   & 0    \\
      a_1/2   & a_0 \\
      a_2/2   & a_1     & a_0 \\
    \vdots    &         &         & \ddots \\
    a_{n}/2   & a_{n-1} & a_{n-2} &        &  a_{0}
   \end{array}\right].
\end{equation}
Note that the first column of $E_n$ is divided by two.
This allows for the immediate computation
$$
   \mbox{Tr}(A_n^{-1}\left[c_{i+j}\right]_{0\leq i,j \leq n}) = (2n+1)a_0.
$$
So the third term (\ref{term3}) is 
\begin{equation}
\label{term3simplified}
  4(2x^2-6x+3)\mbox{Tr}(A_n^{-1}[c_{i+j}]_{0 \leq i,j \leq n})
     = 4(2n+1)(2x^2-6x+3)
\end{equation}
The term
(\ref{term4}) that involves
$$
   \left[c_{i+j-1}\right]_{0 \leq i,j \leq n}
$$
case is similar.  We have the equation
\begin{equation}
\label{FA}
   \left[c_{i+j-1}\right]_{0 \leq i,j \leq n}
      =  F_n A_n + A_n F_n^T
\end{equation}
where
\begin{equation}
\label{Fmatrix}
   F_n = F_n (x) =
        \left[\begin{array}{cccccc}
           0   \\
           a_0    & 0 \\
           a_1    & a_0     & 0 \\
             \vdots               &        & \ddots \\
           a_{n-1}& a_{n-2} & \ldots &  a_0 & 0
         \end{array}\right].
\end{equation}
The identity
(\ref{FA})
leads to the computation
$$
   \mbox{Tr}(A_n^{-1}[c_{i+j-1}]_{0 \leq i,j \leq n}) = 0,
$$
qo that term (\ref{term4})
evaluates to zero:
\begin{equation}
\label{term4simplified}
   -27(2x^2-6x+3)\mbox{Tr}(A_n^{-1}[c_{i+j-1}]_{0 \leq i,j \leq n}) = 0
\end{equation}

Adding the  expressions 
(\ref{term1simplified}), 
(\ref{term2simplified}), 
(\ref{term3simplified}), 
(\ref{term4simplified})
and multiplying through by $H_n$
we obtain the following expression for the first derivative
\begin{lemma}
\label{firstder}
\begin{equation}
\label{det01Eq31}
    (x-3)(2x-3)(4x-3)d_x H_n =  \left( 8nx^2-6(5n+1)x-3(9n^2+13n+8) \right) H_n +2(4n+3) K_n
\end{equation}
\end{lemma}

We next state and prove a lemma which is preparatory to the calculation of the 
second derivative of $H_n$.
First
define two new determinants as follows:
\begin{equation}
\label{defMn}
   M_n =M_n (x)  = \det \left[\begin{array}{ccccccc}
       a_0     & a_1     & \ldots & a_{n-2} &a_{n+1}  & a_{n} \\
       a_1     & a_2     & \ldots & a_{n-1} &a_{n+2}    & a_{n+1} \\
       \vdots  &         & \ddots                        \\
       a_{n-1} & a_n     & \ldots & a_{2n-3}& a_{2n} & a_{2n-1} \\
       a_{n} & a_{n+1} & \ldots & a_{2n-2 } & a_{2n+1}   & a_{2n}
              \end{array}\right]
\end{equation}
\begin{equation}
\label{defNn}
   N_n =N_n (x) = \det \left[\begin{array}{ccccccc}
       a_0     & a_1     & \ldots & a_{n-2} & a_{n-1}  & a_{n+2} \\
       a_1     & a_2     & \ldots & a_{n-1} & a_{n}    & a_{n+3} \\
       \vdots  &         & \ddots  &                        \\
       a_{n-1} & a_n     & \ldots & a_{2n-3} & a_{2n-2} & a_{2n+1} \\
       a_{n} & a_{n+1} & \ldots & a_{2n-2} & a_{2n-1}   & a_{2n+2}
              \end{array}\right]
\end{equation}
where $a_k = a_k(x)$ are the polynomials defined in (\ref{akintro31}).
Then there is a linear relationship between the four determinants 
$H_n, K_n, M_n,$ and $N_n$ as stated in the following lemma.

\begin{lemma}
\label{Delta23LinearRelationship}
\begin{eqnarray}
\nonumber
    & & \mbox{\ \ } 4(4n+5)(x-1)N_n  + ~4(4n+1)(x-1)M_n  \\
\nonumber
         & & + \left( -16 (4n+1)x^3 +8 (36n +7) x^2 -108 (5n +2) x  + 6 (54 n + 31) \right) K_n\\
\label{Delta23Linear}
          & & +~
\left( ( 64 n +16) x^4 + (216 n^2 -24n -12) x^3- (972 n^2 +800n+108) x^2 \right. \\
\nonumber
& & \hspace*{2cm} \left. + (1458n^2 + 1770n +378) x - 729 n^2 -1083 n -324 \right) H_n =0
\end{eqnarray}
\end{lemma}
\noindent
{\em Proof of the Lemma: ~~} First we make use of the recursion in 
in Lemma \ref{lemma2} for each index $k$.
Putting these all together in matrix form, 
we apply the operator 
$$
    \mbox{Tr}(A_n^{-1}*)
$$
to obtain the trace identity 
{\small
\begin{eqnarray}
\label{newterm1}
     &&\hspace*{-5mm}  \mbox{Tr}(A_n^{-1}[4(2(i+j)+5)(x-1)a_{i+j+2}]_{0 \leq i,j \leq n} ) \\
\label{newterm2}
                  && + ~\mbox{Tr}(A_n^{-1}
                             [ (- 2(16x^3-72x^2+135x-81)(i+j) 
                                 -2(24x^3-92x^2+180x-117))a_{i+j+1}
                             ]_{0 \leq i,j \leq n}) \\
\label{newterm3}
          & &+~\mbox{Tr}(A_n^{-1}[(27(2x-3)^3(i+j)
                                +54(2x-3)(2x^2-4x+3))a_{i+j}]_{0 \leq i,j \leq n}) \\
\label{newterm4}
           & &+~ \mbox{Tr}(A_n^{-1}[8(x-1)(2x^2-6x+3)c_{i+j+1}]_{0 \leq i,j \leq n}) \\
\label{newterm5}
           & &+~\mbox{Tr}(A_n^{-1}[2(8x^4-114x^3+324x^2-297x+81)c_{i+j}]_{0 \leq i,j \leq n}) \\
\label{newterm6}
           & &+~\mbox{Tr}(A_n^{-1}[-27x(2x-3)(2x^2-12x+9)c_{i+j-1}]_{0 \leq i,j \leq n}) 
            ~ = ~  0
\end{eqnarray}
}

Each of these six traces (\ref{newterm1})-(\ref{newterm6}) is calculated in a similar manner as was done
above in the calculation of the first derivative of $H_n$.  The
first and fourth trace will involve a small extension to what was used
above.  We will start with the computation of (\ref{newterm1}).
As before we write this as
\begin{equation}
\label{4x}
      4(x-1)\mbox{Tr}(A_n^{-1}((2L_n+\frac{5}{2}I)[a_{i+j+2}]_{0\leq i,j \leq n}  +
                              [a_{i+j+2}]_{0 \leq i,j \leq n} (2L_n+\frac{5}{2}I))).
\end{equation}
Using the fact that
$$\begin{array}{l@{}l}
   4(x-1)\mbox{Tr}&(A_n^{-1}(2L_n+\frac{5}{2}I )[a_{i+j+2}]_{0 \leq i,j \leq n}) \\
         & = 4(x-1)\mbox{Tr}((2L_n+\frac{5}{2}I)[a_{i+j+2}]_{0 \leq i,j \leq n} A_n^{-1})
\end{array}$$
and
$$
   \left((2L_n+\frac{5}{2}I)[a_{i+j+2}]_{0 \leq i,j \leq n} A_n^{-1}\right)^T =
          A_n^{-1}[a_{i+j+2}]_{0 \leq i,j \leq n} (2L_n+\frac{5}{2}I)
$$
we see that we can write (\ref{4x}) as 
\begin{equation}
\label{8x}
      8(x-1)\mbox{Tr}(A_n^{-1}[a_{i+j+2}]_{0 \leq i,j \leq n} (2L_n+\frac{5}{2}I)).
\end{equation}
We will obtain a representation of the trace above in terms of 
$H_n, M_n$, and $N_n$.
Introduce $\tau_0$, $\tau_1$, $\ldots$, $\tau_n$ as follows:
$$
\left[\begin{array}{c}
            \tau_0 \\
            \tau_1 \\
            \vdots \\
            \tau_n
                \end{array}\right]
=
    A_n^{-1}\left[\begin{array}{c}
        a_{n+2} \\
        a_{n+3} \\
        \vdots \\
        a_{2n+2}
           \end{array}\right] .
$$
As before we observe that
$$
   A_n^{-1}[a_{i+j+2}]_{ 0 \leq i,j \leq n}
     = \left[\begin{array}{cccccccccc}
  0 &   &        &        &        &    & \sigma_0     & \tau_0      \\
  0 & 0 &        &        &        &    & \sigma_1     & \tau_1      \\
  1 & 0 & \ddots &        &        &    & \sigma_2     & \tau_2      \\
    & 1 &        & \ddots &        &    & \sigma_3     & \tau_3      \\
    &   & \ddots &        & \ddots &    & \vdots       & \vdots      \\
    &   &        & \ddots &        &  0 & \sigma_{n-2} & \tau_{n-2}  \\
    &   &        &        & \ddots &  0 & \sigma_{n-1} & \tau_{n-1}  \\
    &   &        &        &        &  1 & \sigma_n     & \tau_n
             \end{array}\right]
$$
and therefore
$$
    A_n^{-1}[a_{i+j+2}]_{ 0 \leq i,j \leq n} (2L_n+\frac{5}{2}I)
$$
expands to
$$
   \left[\begin{array}{cccccccccc}
  0  &        &        &        &         & (4n+1)\sigma_0/2     & (4n+5)\tau_0/2 \\
  0  & \ddots &        &        &         & (4n+1)\sigma_1/2     & (4n+5)\tau_1/2      \\
 5/2 &        & \ddots &        &         & (4n+1)\sigma_2 /2    & (4n+5)\tau_2/2      \\
     & \ddots &        & \ddots &         & \vdots               & \vdots              \\
     &        & \ddots &        &    0    & (4n+1)\sigma_{n-2}/2 & (4n+5)\tau_{n-2}/2  \\
     &        &        & \ddots &    0    & (4n+1)\sigma_{n-1}/2 & (4n+5)\tau_{n-1}/2  \\
     &        &        &        & (4n-3)/2& (4n+1)\sigma_n/2     & (4n+5)\tau_n/2
             \end{array}\right].
$$
Thus (\ref{8x}) 
simplifies to
$$
        4(x-1)((4n+1)\sigma_{n-1} + (4n+5)\tau_n)
$$
By Cramer's rule we have
\begin{eqnarray}
\label{Cramersrule2}
  \sigma_{n-1} & = & \frac{M_n}{H_n} \\
\label{Cramersrule3}
  \tau_n       & = & \frac{N_n}{H_n}\;
\end{eqnarray}
so, in summary, the first trace (\ref{newterm1}) evaluates to
\begin{equation}
\label{newterm1simplified}
         4(x-1)((4n+1) \frac{M_n}{H_n} + (4n+5) \frac{N_n}{H_n}).
\end{equation}
Now we consider the second trace (\ref{newterm2}).
The evaluation of this trace proceeds exactly as the evaluation of
the trace (\ref{term1}) in the computation of the derivative of $H_n$.  (\ref{newterm2}) 
evaluates to 
\begin{equation}
\label{newterm2simplified}
   - 2 \left( 2(16x^3-72x^2+135x-81)n+(24x^3-92x^2+180x-117) \right) \frac{K_n}{H_n}.
\end{equation}
Similarly, the third trace (\ref{newterm3}) 
evaluates to
\begin{equation}
\label{newterm3simplified}
    (n+1) \left( 27(2x-3)^3n +54(2x-3)(2x^2-4x+3) \right) .
\end{equation}

The next trace (\ref{newterm4})
requires that we define a matrix $G_n$ in a similar manner to $E_n$
and $F_n$.  Specifically we define
$$
   G_n = G_n(x) = \left[\begin{array}{ccccccccc}
                  a_1/2  & a_0/2 \\
                  a_2/2  & a_1/2       & a_0 \\
                  a_3/2  & a_2/2       & a_1    & a_0 \\
 \vdots                         &    \vdots          &        &        & \ddots \\
              a_{n}/2    & a_{n-1}/2   & a_{n-2}& \ldots &        & a_0 \\
              a_{n+1}/2  & a_n/2       & a_{n-1}& \ldots &        & a_1
               \end{array}\right].
$$
With this definition of $G_n$ we can write
$$
  [c_{i+j+1}]_{0 \leq i,j \leq n}  =  G_n A_n + A_n G_n^T 
   + a_0 
    \left[\begin{array}{ccccc}
          0&  &  \ldots & 0 \\
       \vdots &     &        & \vdots \\
          0 & 0 &  \ldots & 0 \\
          a_{n+1} & a_{n+2}    &  \ldots       & a_{2n+1}
       \end{array}\right] 
   + a_0 
    \left[\begin{array}{ccccc}
          0&    \ldots & 0&  a_{n+1}\\
           \vdots &    &  0& a_{n+2}\\
       &  &      & \vdots \\
          0 &  \ldots  &0     & a_{2n+1}
       \end{array}\right] 
$$
We have
$$
   A_n^{-1}
    \left[\begin{array}{ccccc}
          0&    \ldots & 0&  a_{n+1}\\
           \vdots &    &  0& a_{n+2}\\
       &  &      & \vdots \\
          0 &  \ldots  &0     & a_{2n+1}
       \end{array}\right] 
     = 
\left[\begin{array}{ccccc}
          0&    \ldots & 0&  \sigma_0 \\
           \vdots &    &  0& \sigma_1 \\
       &  &      & \vdots \\
          0 &  \ldots  &0     & \sigma_n
       \end{array}\right]
$$
where $ \sigma_i$ is as defined in (\ref{sigman}). It follows that
$$
  \mbox{Tr}\left(
      A_n^{-1}
 \left[\begin{array}{ccccc}
          0&    \ldots & 0&  a_{n+1}\\
           \vdots &    &  0& a_{n+2}\\
       &  &      & \vdots \\
          0 &  \ldots  &0     & a_{2n+1}
       \end{array}\right]
\right) = \sigma_{n}.
$$
The term
$$
  \mbox{Tr}\left(
      A_n^{-1}
\left[\begin{array}{ccccc}
          0&  &  \ldots & 0 \\
       \vdots &     &        & \vdots \\
          0 & 0 &  \ldots & 0 \\
          a_{n+1} & a_{n+2}    &  \ldots       & a_{2n+1}
       \end{array}\right]
\right) $$
comes out to $\sigma_n$ also because
$$
  \left(
\left[\begin{array}{ccccc}
          0&  &  \ldots & 0 \\
       \vdots &     &        & \vdots \\
          0 & 0 &  \ldots & 0 \\
          a_{n+1} & a_{n+2}    &  \ldots       & a_{2n+1}
       \end{array}\right]
A_n^{-1}\right)^T
     =    A_n^{-1}
 \left[\begin{array}{ccccc}
          0&    \ldots & 0&  a_{n+1}\\
           \vdots &    &  0& a_{n+2}\\
       &  &      & \vdots \\
          0 &  \ldots  &0     & a_{2n+1}
       \end{array}\right]
$$
Putting this together we see that the fourth trace (\ref{newterm4})
evaluates to
\begin{equation}
\label{newterm4simplified}
  16(x-1)(2x^2-6x+3) \left( n(x+4) + \frac{K_n}{H_n} \right) .
\end{equation}

The evaluation of the fifth and sixth trace (\ref{newterm5}) and (\ref{newterm6}) 
are done in exactly the same manner as the traces (\ref{term3}) and (\ref{term4}).
For
the fifth trace (\ref{newterm5})
we obtain
\begin{equation}
\label{newterm5simplified}
 \begin{array}{l@{}l}
    \mbox{Tr}(A_n^{-1}[2(8x^4-114x^3+&324x^2-297x+81)c_{i+j}]_{0 \leq i,j \leq n}) \\
              & = 2(2n+1)(8x^4-114x^3+324x^2-297x+81)
 \end{array}
\end{equation}
and
the sixth trace (\ref{newterm6}) evaluates to zero just as the trace in (\ref{term4}):
\begin{equation}
\label{newterm6simplified}
   \mbox{Tr}(A_n^{-1}[- 27x(2x-3)(2x^2-12x+9)c_{i+j-1}]_{0 \leq i,j \leq n}) = 0 ~.
\end{equation}
Adding the expression we have found  for the six traces
in (\ref{newterm1simplified})-(\ref{newterm6simplified}) we obtain the identity in
(\ref{Delta23Linear}).

\subsection{Calculating the Second Derivative}
\label{seconddersection}

We are now ready to calculate the second derivative of $H_n$.
We begin with 
equation (\ref{det01Eq31}) for the derivative.  The
first step is to
replace $K_n$ in equation (\ref{det01Eq31}) with the
representation of $K_n$ as a trace from (\ref{DeltaOneAsTrace31}).
Inserting this
expression into equation (\ref{det01Eq31}) we get
\begin{equation}
\label{begin2der}
    (x-3)(2x-3)(4x-3)d_x H_n =  (8nx^2-6(5n+1)x-3(9n^2+13n+8)) H_n + 2(4n+3)\mbox{Tr}(A_n^{-1}B_n)H_n
\end{equation}
In order to obtain the second derivative of $H_n$, we 
differentiate (\ref{begin2der}). We obtain
\begin{eqnarray*}
    (x-3)(2x-3)(4x-3)  d_x^2 H_n  + &&\\
         (8(n+3)x^2-6(5n-13)x-3(9n^2+13n+29))  d_x H_n + &&\\
         2(8n x-15n-3)  H_n ~~~~~ && \\
       = ~~2(4n+3)(d_x\mbox{Tr}(A_n^{-1}B_n)) H_n + 2(4n+3)\mbox{Tr}(A_n^{-1}B_n)d_x H_n && . 
\end{eqnarray*}
Multiply both sides of the equation by
$$
   (x-3)(2x-3)(4x-3)
$$
The second trace term on the right hand side now becomes  
\begin{equation}
\label{sectrace}
     2(4n+3)(x-3)(2x-3)(4x-3)\mbox{Tr}(A_n^{-1}B_n)d_x H_n
\end{equation}
Using the identities 
(\ref{det01Eq31}) and (\ref{DeltaOneAsTrace31}),
(\ref{sectrace}) can be written 
in terms of $H_n$ and $K_n$ as follows:
as
$$
     2(4n+3) \left((8nx^2-6(5n+1)x-3(9n^2+13n+8)) K_n
                + 2(4n+3) \frac{K_n^2}{H_n}\right).
$$
Substituting back, we get
\begin{eqnarray}\nonumber
  (x-3)^2(2x-3)^2(4x-3)^2d_x^2H_n - && \\ \nonumber
(x-3)(2x-3)(4x-3) (8(n-3)x^2-6(5n-13)x -3(9n^2+13n+29))d_x H_n -  && \\ \nonumber
         2(x-3)(2x-3)(4x-3)(8n x-15n-3) H_n - && \\
\label{ExpandedSecondDerivativeNoSimplify31}
          2(4n+3)(8nx^2-6(5n+1)x -3(9n^2+13n+8))K_n - && \\ \nonumber
         4(4n+3)^2\frac{K_n^2}{H_n} - &&  \\ \nonumber
   2(4n+3)(x-3)(2x-3)(4x-3) \left( d_x\mbox{Tr}(A_n^{-1}B_n) \right) H_n ~=~0&&
\end{eqnarray}
We will now focus on simplification of the derivative of the trace 
$$
    (x-3)(2x-3)(4x-3)(d_x\mbox{Tr}(A_n^{-1}B_n))
$$
which is a factor of the last term on the left of equation
(\ref{ExpandedSecondDerivativeNoSimplify31}).
We use the fact that
$$
    d_x A_n^{-1} = - A_n^{-1}(d_x A_n)A_n^{-1}
$$
and write
\begin{equation}
\label{trderivative}
    d_x\mbox{Tr}(A_n^{-1}B_n)   = 
    \mbox{Tr}(A_n^{-1}d_x B_n) - 
\mbox{Tr}(A_n^{-1}(d_xA_n)A_n^{-1}B_n)
\end{equation}

Using equation (\ref{FirstId31}) from Lemma \ref{lemma1},
and multiplying equation (\ref{trderivative}) by $(x-3)(2x-3)(4x-3)$,
we can write
$$
 (x-3)(2x-3)(4x-3)\mbox{Tr}(A_n^{-1}d_x B_n)
$$
as
\begin{eqnarray}
\label{tr1term1}
 && \mbox{\ \ \ } \mbox{Tr}(A_n^{-1}[2(2(i+j+1)+3)a_{i+j+2}]_{0 \leq i,j \leq n}) \\
\label{tr1term2}
 && \mbox{}+\mbox{Tr}(A_n^{-1}[- (8x^2-18x+27(i+j+1)+36)a_{i+j+1}]_{0 \leq i,j \leq n}) \\
\label{tr1term3}
 && \mbox{}+\mbox{Tr}(A_n^{-1}[4(2x^2-6x+3)c_{i+j+1}]_{0\leq i,j \leq n}) \\
\label{tr1term4}
 && \mbox{}+\mbox{Tr}(A_n^{-1}[-27(2x^2-6x+3)c_{i+j}]_{0 \leq i,j \leq n}) 
\end{eqnarray}
and we can write
$$
 (x-3)(2x-3)(4x-3) \mbox{Tr}(A_n^{-1}(d_xA_n)A_n^{-1}B_n)
$$
as
\begin{eqnarray}
\label{tr2term1}
 &&\mbox{\ \ \ } \mbox{Tr}(A_n^{-1}[2(2(i+j)+3)a_{i+j+1}]_{0 \leq i,j \leq n} A_n^{-1}B_n)    \\
\label{tr2term2}
 && \mbox{}+\mbox{Tr}(A_n^{-1}[- (8x^2-18x+27(i+j)+36)a_{i+j}]_{0 \leq i,j \leq n} A_n^{-1}B_n) \\
\label{tr2term3}
 && \mbox{}+\mbox{Tr}(A_n^{-1}[4(2x^2-6x+3)c_{i+j}]_{0 \leq i,j \leq n} A_n^{-1}B_n)    \\
\label{tr2term4}
 && \mbox{}+\mbox{Tr}(A_n^{-1}[-27(2x^2-6x+3)c_{i+j-1}]_{0 \leq i,j \leq n} A_n^{-1}B_n)
\end{eqnarray}
We have already evaluated traces that are similar to the first four terms 
(\ref{tr1term1})-(\ref{tr1term4})
in the calculation of the trace terms in (\ref{term1})-(\ref{term4}). In this case we obtain
{\small
\begin{eqnarray}
\label{tr1term1simplified}
     \mbox{Tr}(A_n^{-1}[2(2(i+j+1)+3)a_{i+j+2}]_{0 \leq i,j \leq n}) & = &
              2((4n+1) \frac{M_n}{H_n} + (4n+5) \frac{N_n}{H_n})\\
\label{tr1term2simplified}
  \mbox{Tr}(A_n^{-1}[-(8x^2-18x+27(i+j+1)+36)a_{i+j+1}]_{0 \leq i,j \leq n}) 
                              &=& -(8x^2-18x+54n+63)\frac{K_n}{H_n}\\
\label{tr1term3simplified}
   \mbox{Tr}(A_n^{-1}[4(2x^2-6x+3)c_{i+j+1}]_{0 \leq i,j \leq n}) &=&
              (2x^2-6x+3)(8n(x+4) + 8 \frac{K_n}{H_n})\\
\label{tr1term4simplified}
   \mbox{Tr}(A_n^{-1}[-27(2x^2-6x+3)c_{i+j}]_{0 \leq i,j \leq n}) & =& -27(2n+1)(2x^2-6x+3)
\end{eqnarray}
}
Next we evaluate the traces (\ref{tr2term1})-(\ref{tr2term4}). 
To calculate the first of these (\ref{tr2term1}) 
We expand the matrix
$$
   A_n^{-1}[2(2(i+j)+3)a_{i+j+1}]_{0 \leq i,j \leq n} A_n^{-1}B_n
$$
as
$$
   A_n^{-1}\left((4L_n+3I)B_n+B_n(4L_n+3I)\right)A_n^{-1}B_n~.
$$
The trace (\ref{tr2term1}) thus can be evaluated as 
$$\begin{array}{r@{}l}
   \mbox{Tr}((4L_n+3I)B_n A_n^{-1}B_n A_n^{-1}) +& \mbox{Tr}(A_n^{-1}B_n (4L_n+3I)A_n^{-1}B_n) \\
               = & 2 \mbox{Tr}(A_n^{-1}B_n (4L_n+3I)A_n^{-1}B_n) \\
               = & 2 \mbox{Tr}((4L_n+3I)\left(A_n^{-1}B_n \right)^2)
\end{array}
$$
where we used the symmetry of the matrices in the trace calculation.
Using the expression (\ref{A1AsSigma}) for $ A_n^{-1}B_n$, we have
$$
   \left(A_n^{-1}B_n \right)^2  = 
   \left[\begin{array}{cccccccccccc}
 0 &        &         &        &   & \sigma_0     & \sigma_0\sigma_n \\
 0 & \ddots &         &        &   & \sigma_1     & \sigma_0+\sigma_1\sigma_n \\
 1 &        & \ddots  &        &   & \sigma_2     & \sigma_1+\sigma_2\sigma_n \\
   &        &         & \ddots &   & \vdots       &   \vdots         \\
   &        &         &        & 0 & \sigma_{n-2} & \sigma_{n-3}+\sigma_{n-2}\sigma_n \\
   &        &         &        & 0 & \sigma_{n-1} & \sigma_{n-2}+\sigma_{n-1}\sigma_n \\
   &        &         &        & 1 & \sigma_{n}   & \sigma_{n-1}+\sigma_n^2 \\
        \end{array}\right].
$$
Therefore
$$
    (4L_n+3I)\left(A_n^{-1}B_n \right)^2
       = \left[\begin{array}{cccccccccccc}
 0 &        &         &        &   & 3\sigma_0     & 3\sigma_0\sigma_n \\
 0 & \ddots &         &        &   & 7\sigma_1     & 7(\sigma_0+\sigma_1\sigma_n) \\
 1 &        & \ddots  &        &   & 10\sigma_2    & 10(\sigma_1+\sigma_2\sigma_n) \\
   &        &         & \ddots &   & \vdots       &   \vdots         \\
   &        &         &        & 0 & (4n-5)\sigma_{n-2} 
                                                  & (4n-5)(\sigma_{n-3}+\sigma_{n-2}\sigma_n) \\
   &        &         &        & 0 & (4n-1)\sigma_{n-1} 
                                                  &(4n-1) (\sigma_{n-2}+\sigma_{n-1}\sigma_n) \\
   &        &         &        & 1 & (4n+3)\sigma_{n} 
                                                  &(4n+3)(\sigma_{n-1}+\sigma_n^2) \\
        \end{array}\right].
$$
and thus
\begin{eqnarray}
\nonumber
      \mbox{Tr}(A_n^{-1}[2(2(i+j)+3)a_{i+j+1}]_{\leq i,j \leq n} A_n^{-1}B_n)
    &= & 2 ( (4n-1)\sigma_{n-1}+(4n+3)(\sigma_{n-1}+\sigma_n^2))\\
\label{tr2term1simplified}
       & =& 4(4n+1)\frac{M_n}{H_n}+2(4n+3)\frac{K_n^2}{H_n^2}
\end{eqnarray}
where the last equality is a consequence of the identities 
(\ref{Cramersrule1}) and (\ref{Cramersrule2}).

Using similar reasoning, we see that the trace term (\ref{tr2term2})
$$
  \mbox{Tr}(A_n^{-1}[-(8x^2-18x+27(i+j)+36)a_{i+j}]_{0 \leq i,j \leq n} A_n^{-1}B_n),
$$
evaluates to
\begin{equation}
\label{tr2term2simplified}
    -(8x^2-18x+54n+36)\frac{K_n}{H_n}.
\end{equation}

The trace term (\ref{tr2term3})
is expanded using (\ref{M0Equation}).  We
obtain the following trace identities
\begin{eqnarray*}
  \mbox{Tr}(A_n^{-1}[4(2x^2-6x+3)c_{i+j}]_{0 \leq i,j \leq n} A_n^{-1}B_n) &=&
     4(2x^2-6x+3)\mbox{Tr}(A_n^{-1}(E_nA_n + A_nE_n^T)A_n^{-1}B_n) \\
    & = & 4(2x^2-6x+3)(\mbox{Tr}(A_n^{-1}E_nB_n) 
                          +\mbox{Tr}(E_n^TA_n^{-1}B_n)) \\
    & = & 8(2x^2-6x+3) \mbox{Tr}(E_n^TA_n^{-1}B_n).
\end{eqnarray*}
Since
$$
  E_n^TA_n^{-1}B_n = \left[\begin{array}{cccccc}
   a_1/2 & a_2/2  & \ldots & a_n/2     & a_0\sigma_0 /2 +\ldots \\
    a_0  &  a_1   & \ldots & a_{n-1}   & a_0\sigma_1+\ldots \\
         & \ddots &        &           & \vdots \\
         &        & \ddots &    a_1       & a_0 \sigma_{n-1}+a_1 \sigma_n \\
         &  0      &        & a_0       & a_0\sigma_n
                   \end{array}\right]
$$
the trace reduces as follows:
\begin{eqnarray} \nonumber
  \mbox{Tr}(A_n^{-1}[4(2x^2-6x+3)c_{i+j}]_{0 \leq i,j \leq n} A_n^{-1}B_n)
    & = & 4(2x^2-6x+3)((2n-1)a_1 + 2 a_0\sigma_n) \\ \nonumber
    & = & 4(2n-1)(x+4)(2x^2-6x+3) \\
& & \hspace*{1cm} + 8(2x^2-6x+3)\frac{K_n}{H_n}.
\label{tr2term3simplified}
\end{eqnarray}

Using the expansion (\ref{FA}), 
the final trace term (\ref{tr2term4})
simplifies as follows:
\begin{equation}
\label{tr2term4simplified}
  \mbox{Tr}(A_n^{-1}[-27(2x^2-6x+3)c_{i+j-1}]_{0\leq i,j\leq n} A_n^{-1}B_n)
      = - 54n (2x^2-6x+3) 
\end{equation}

Now we return to equation 
(\ref{ExpandedSecondDerivativeNoSimplify31}). In this equation, there is a term 
containing
$$
   d_x\mbox{Tr}(A_n^{-1}B_n)
$$
We expanded the derivative above 
as the difference of two traces in 
(\ref{trderivative}). We multiplied (\ref{trderivative}) by
$ (x-3) (2x-3)(4x-3)$ and expressed the product
$$ 
(x-3) (2x-3)(4x-3) d_x\mbox{Tr}(A_n^{-1}B_n))
$$
as the sum of four terms in
(\ref{tr1term1}),
(\ref{tr1term2}),
(\ref{tr1term3}),
(\ref{tr1term4})
minus the sum of four terms in
(\ref{tr2term1}),
(\ref{tr2term2}),
(\ref{tr2term3}),
(\ref{tr2term4}). 
These eight traces in turn
were simplified as 
(\ref{tr1term1simplified}),
(\ref{tr1term2simplified}),
(\ref{tr1term3simplified}),
(\ref{tr1term4simplified}),
(\ref{tr2term1simplified}),
(\ref{tr2term2simplified}),
(\ref{tr2term3simplified}),
(\ref{tr2term4simplified}). 
Computing the sum of 
(\ref{tr1term1simplified}),
(\ref{tr1term2simplified}),
(\ref{tr1term3simplified}),
(\ref{tr1term4simplified})
minus the sum of 
(\ref{tr2term1simplified}),
(\ref{tr2term2simplified}),
(\ref{tr2term3simplified}),
(\ref{tr2term4simplified})
we obtain
\begin{eqnarray*}
H_n (x-3) (2x-3)(4x-3) d_x\mbox{Tr}(A_n^{-1}B_n)) &= & 
2(4 n + 5) N_n  - 2(4 n + 1) M_n  - 27  K_n +  \\
& &
      (4 x  - 11 )(2 x^2  - 6 x + 3 )H_n - 2( 4 n+3) \frac{K_n^2}{H_n}  
\end{eqnarray*}
Now we multiply this through by $ 2 (4n+3)$. 
The left hand side becomes the last term on the 
left hand side of equation
(\ref{ExpandedSecondDerivativeNoSimplify31}).
Substituting 
and simplifying, we obtain
the intermediate identity
\begin{eqnarray}
\nonumber
 (x-3)^2(2x-3)^2(4x-3)^2d_x^2 H_n - && \\
\nonumber
(x-3)(2x-3)(4x-3)(8(n-3)x^2-6(5n-13)x - 3( 9 n^2 +13 n+29))d_x H_n - && \\
\nonumber
2( 64n x^4 -424n x^3 +2 (475 n-6)x^2 -3 (283 n-15)x + 3 (91 n-6)) H_n - &&\\
\label{intermediateidentity}
2(4n+3)(8n x^2 -6(5n+1) x-3 (9n^2 +13n +17)) K_n + &&\\
\nonumber
 4(4n+1)(4n+3) M_n - &&\\
\nonumber
 4(4n+3)(4n+5)N_n ~= ~ 0&&.
\end{eqnarray}

Recall that equation (\ref{det01Eq31}) of Lemma  \ref{firstder}
expresses 
$d_x H_n$ in terms of $H_n$ and $K_n$.  Also,
Lemma \ref{Delta23LinearRelationship} gives a linear relationship
between $H_n$, $K_n$, $M_n$ and $N_n$.  Thus by
using these two equations we can eliminate both the $d_x H_n$ and the
$M_n$ term in the above formula:

\begin{eqnarray}
\nonumber
4 (4n+1)(x-1)(x-3)^2(2x-3)^2(4x-3)^2 d_x^2 H_n  -&& \\
\nonumber
4 (4n+1) (64 n (n-1) x^5 -96(3 n^2 -8n -2) x^4 +12 (36n^3 +163 n^2 -65 n-8) x^3 -  &&\\
\nonumber
4 (459 n^3 +1661 n^2 +949 n +417)x^2 + 
  3 (243 n^4 +2106 n^3 +5192n^2 +4291 n +1482) x && \\
\nonumber
-3 (243 n^4+ 1674n^3+3679 n^2 +3140 n+1008)) H_n  +&&\\
\label{SecondDerivativeAs013}
8 (4n+1) (4n+3) (16(n+2)x^3 -4 (17n+31) x^2 + 6(9n^2+48n+53) x - && \\
\nonumber 
3 (18n^2+80n+77) ) K_n  - &&\\
\nonumber
32(4n+1)(4n+3)(4n+5)(x-1) N_n ~=~ 0&&.
\end{eqnarray}

At this point we have enough information on the interrelationship between 
$ H_n, K_n, N_n$ and $M_n$ that allows us to
derive a number of evaluations of $H_n(x)$ for special values of $x$. These evaluations also turn out to 
be essential for the final step of the proof of Theorem \ref{thm1} in Section \ref{lastpiece}.

\section{Product form evaluations of $H_n(x)$ at special $x$}
\label{evaluatex}

Our calculations for the special values rely on identities we have proved for
$ H_n, K_n, N_n$ and $M_n$ 
as well as a
Dodgson-like determinant identity which is useful.
This identity is as follows:
\begin{proposition}
\label{detidentprop}
Suppose the determinants 
$H_n, K_n, M_n$, and $N_n$ are defined as in
(\ref{defHn}), 
(\ref{defKn}),
(\ref{defMn}), and 
(\ref{defNn}) respectively.
Then
\begin{equation}
\label{detident} 
  H_{n-1} H_{n+1} = H_n N_n - H_n M_n - K_n^2.
\end{equation}
\end{proposition}
\begin{proof}
(Proof of Proposition \ref{detidentprop})
This identity can be proved using techniques
similar to those in \cite{Knuth96}.
We first prove a general determinant expansion result, and then specialize to Hankel determinants to
obtain (\ref{detident}).
Consider two matrices 
\begin{eqnarray*}
    R_n&  =&  [r_{i,j}]_{0\le i,j\le n} \\
   X_{n+1}&  =&  [x_{i,j}]_{0\le i,j\le n+1}
\end{eqnarray*}
where ultimately we will set for all $ i, j$
$$
   r_{i,j} = x_{i,j} = a_{i+j}(x).
$$
Consider the sum
\begin{equation}\label{detsum}
  \sum_{k=0}^{n+1} (-1)^{n+1-k}
\left(
\det \left[\left\{\begin{array}{ll}r_{i,j}   & \mbox{if $i\neq n$} \\
                                      x_{k,j+1} & \mbox{if $i = n$}
                     \end{array} \right.
\right]_{0\le i,j\le n} 
     \det 
        \left[\left\{\begin{array}{ll}x_{i,j}   & \mbox{if $i < k$} \\
                                      x_{i+1,j}   & \mbox{if $i \ge k$}
                     \end{array}\right]_{0\le i,j\le n} \right.
\right)
\end{equation}
as a function of the matrix $X = X_{n+1}$.
It is not difficult to prove that if two adjacent rows of $X$ are
switched then the sum (\ref{detsum}) changes sign. Since
the set of all permutations of the
rows is generated by flips of adjacent rows, it follows that (\ref{detsum}) is alternating.

In addition to this (\ref{detsum}) is linear in each row of $X$.
Since (\ref{detsum}) is both alternating and multilinear, it is equal to 
a multiple (depending on $R_n$) of $\det( X_{n+1})$.
To determine the multiple we set $X$ to the matrix:
$$
   X_{n+1}  = \left[
      \begin{array}{cccccc}
   1 &    0    & \ldots  &        &   0     \\
   0 & r_{0,0} & r_{0,1} & \ldots & r_{0,n} \\
   0 & r_{1,0} & r_{1,1} & \ldots & r_{1,n} \\
 \vdots                                     \\
   0 & r_{n,0} & r_{n,1} & \ldots & r_{n,n} \\
      \end{array}\right]
$$
In this case 
$$
  \det ( X_{n+1} ) = \det ( R_n) ~.
$$
In the sum~(\ref{detsum}) only the term corresponding to $ i = n+1$ survives and this term 
itself evaluates to
$$
  \det(R_n) \det (R_{n-1}).
$$
Therefore for all $R$ and $X$, the sum~(\ref{detsum}) evaluates to 
$ \det(R_{n-1}) \det(X_{n+1}).  $
When 
$R$ and $X$ are Hankel, only the last three terms of in the sum (\ref{detsum}) survive. 
In particular, if we set
$$
   r_{i,j} = x_{i,j} = a_{i+j}(x)
$$
we obtain identity (\ref{detident}).
\end{proof}

\subsection{The evaluation of $H_n(3)$}
\label{3section}

We use  the identities 
involving $H_n, K_n, M_n$, and $N_n$ that we already have, 
along with (\ref{detident}).  At the point $x=3$
the left hand side of (\ref{detident}) evaluates to 
$H_{n-1} ( 3) H_{n+1} (3)$.
The right hand side of (\ref{detident}) is more difficult to evaluate; we have to use the identities
from the previous section. At $x=3$, equation (\ref{det01Eq31}) simplifies
to
\begin{equation}
\label{x31}
   2(4n+3)K_n = 3(9n^2+19n+14) H_n.
\end{equation}
Therefore
\begin{equation}
\label{x32}
   K_n (3) = \frac{3(9n^2+19n+14)} {2(4n+3)} H_n(3)
\end{equation}
Now we can use this together with equation (\ref{SecondDerivativeAs013})
to calculate $N_n$ at $x=3$
\begin{equation}
\label{x33}
  N_n (3) = \frac{3(243n^4+1512n^3+3605n^2+4144n+1920)} {8(4n+3)(4n+5)} H_n (3)
\end{equation}
Finally we use the linear relationship of Lemma \ref{Delta23LinearRelationship} to calculate
$M_n$ at $x=3$,
\begin{equation}
\label{x34}
M_n (3) = - \frac{3 n (243 n^3 +540 n^2 +559 n +58)}{8 (4n+1) (4n+3)} H_n (3)
\end{equation}
We can now write the right hand side of (\ref{detident}) as
$$
\frac{ 9 (3n+4)(3n+5)(6n-1)(6n+1)} {4 (4n+1)(4n+3)^2 (4n+5)} {H_n(3)}^2
$$
Letting $H_n = H_n(3)$, we have
\begin{equation}
\label{Crecursion}
H_{n-1} H_{n+1} = 
\frac{ 9 (3n+4)(3n+5)(6n-1)(6n+1)} {4 (4n+1)(4n+3)^2 (4n+5)} H_n^2
\end{equation}

Writing (\ref{Crecursion}) as 
a recursion for $ \frac{H_{n+1}}{H_n}$, we find that 
$$
 \frac{H_{n+1}}{H_n} = \frac{H_1}{H_0} \prod_{i=1}^n \frac{ 9 (3i+4) (3i+5) (6i-1) (6i+1)}{4 (4i+1) (4i+3)^2
(4i+5)}
$$
Thus
$$
H_n=  H_0  (\frac{H_1}{H_0})^{n} \prod_{j=1}^{n-1} \prod_{i=1}^{n-j}\frac{ 9 (3i+4) (3i+5) (6i-1) (6i+1)}{4 (4i+1) (4i+3)^2 (4i+5)}
$$

We compute directly that $H_0 = 1$ and $ H_1 = -1$.
Rearranging
\begin{equation}
\label{Cdoubleproduct}
H_n = (-1)^n \prod_{j=1}^{n-1} \prod_{i=1}^{n-j}\frac{ 9 (3i+4) (3i+5) (6i-1) (6i+1)}{4 (4i+1) (4i+3)^2 (4i+5)}
\end{equation}
This double product expression for $H_n$ can be rewritten in terms of factorials. We record this as a corollary.
\begin{corollary}
\label{3corollary}
Suppose $a_k(x)$ is defined as in 
(\ref{akintro31})
and $H_n(x)= \det [ a_{i+j}(x) ]_{0 \leq i,j \leq n}$. Then
\begin{equation}
\label{Cfactorials}
  H_n(3) = (-1)^n \prod_{i=1}^n \frac{ (6i-3)!(3i+2)!(2i-1)!} {(4i-1)!(4i+1)!(3i-2)!}
\end{equation}
\end{corollary}

\subsection{The evaluation of $H_n(\threehalves)$}
\label{3over2section}

For $ x=\threehalves $, the left hand side of identity 
(\ref{det01Eq31}) of Lemma  \ref{firstder}
vanishes, and we obtain
\begin{equation}
\label{3over21}
K_n (\threehalves ) = \frac{3( 9 n^2 + 22n + 11) } {2( 4 n + 3) } H_n (\threehalves )
\end{equation}
Next, we evaluate 
(\ref{SecondDerivativeAs013}) at $ x =\threehalves $ using this expression for $K_n$. The 
second derivative terms drops out and we obtain
\begin{equation}
\label{3over22}
N_n( \threehalves ) = \frac{3( 243 n^4 + 1674 n^3 + 4031 n^2 +3934 n + 1290 ) }{8( 4 n + 3)( 4n + 5) }
H_n( \threehalves ) 
\end{equation}
Finally, using these expressions for $K_n (\threehalves )$ and $N_n( \threehalves )$ in the linear identity
(\ref{Delta23Linear}) of Lemma \ref{Delta23LinearRelationship} at $ x= \threehalves $, we obtain
\begin{equation}
\label{3over23}
M_n( \threehalves ) = 
\frac{-3 n ( 243 n^3 +702 n^2 +547 n + 118) }{8( 4n +  1) ( 4n +3 ) }
H_n( \threehalves ) 
\end{equation}
Substituting these in (\ref{detident}) with $ H_n = H_n( \threehalves )$, we obtain 
\begin{equation}
\label{3over24}
H_{n-1} H_{n+1} =  \frac{9( 3n + 2) ( 3n + 4) ( 6n + 1) ( 6n + 5)}{4( 4n + 1) {(4n + 3)}^2 (4n + 5)} H_n^2 
\end{equation}
This can be written as a recursion for $ \frac{H_n}{H_{n-1}}$ giving
\begin{equation}
\label{3over25}
H_n = H_{n-1} \cdot \frac{H_1}{H_0} \prod_{i=1}^{n-1} 
  \frac{9( 3i + 2) ( 3i + 4) ( 6i + 1) ( 6i + 5)}{4( 4i + 1) {(4i + 3)}^2 (4i + 5)} 
\end{equation}
Iterating (\ref{3over25}), and using the fact that 
$H_0( \threehalves ) = 1$ and  
$H_1( \threehalves ) = 2$, 
we obtain the product formula
\begin{equation}
\label{3over26}
H_n( \threehalves ) 
= 2^n \cdot  \prod_{j=1}^{n-1} \prod_{i=1}^{n-j} \frac{9( 3i + 2) ( 3i + 4) ( 6i + 1) ( 6i + 5)}{4( 4i + 1) {(4i + 3)}^2 (4i + 5)} 
\end{equation}
The expression (\ref{3over26}) can equivalently be written in terms of factorials as in (\ref{3over27}).
\begin{corollary}
\label{3over2corollary}
Suppose $a_k(x)$ is defined as in 
(\ref{akintro31})
and $H_n(x)= \det [ a_{i+j}(x) ]_{0 \leq i,j \leq n}$. Then
\begin{equation}
\label{3over27}
H_n( \threehalves ) = \prod_{i=1}^n \frac{(2i-1)! (6i)! (3i+1)}{2 (4i-1)! (4i+1)!}
\end{equation}
\end{corollary}

\subsection{The evaluation of $H_n(\threequarters)$}
\label{3over4section}

We can mimic these steps for the evaluation of 
$H_n$ at $ x =\threequarters$ as well. 
Evaluating the 
identities 
(\ref{det01Eq31}) of Lemma  \ref{firstder},
(\ref{SecondDerivativeAs013}), and 
(\ref{Delta23Linear}) of Lemma \ref{Delta23LinearRelationship} at  $ x =\threequarters$, we obtain 
respectively,
\begin{eqnarray}
\label{3over41}
K_n (\threequarters )&=& \frac{3(18 n^2 +38n+  19 ) }
  {4( 4n + 3 ) } H_n (\threequarters )\\
\label{3over42}
N_n (\threequarters )&=& \frac{3( 486 n^4 + 3024 n^3 + 6724 n^2 +6290 n + 2085) }{16( 4n + 3 ) ( 4n+ 5) }H_n
(\threequarters )\\
\nonumber
M_n (\threequarters )&=& \frac{-3n( 243 n^3 + 540n^2 + 316n + 31) }{8( 4n + 1)(4n +  3) } H_n
(\threequarters )
\end{eqnarray}
Substituting these in (\ref{detident}) with $ H_n = H_n( \threequarters )$, 
we obtain the recursion
\begin{equation}
\label{3over43}
H_{n-1} H_{n+1} =  \frac{9( 3n + 1) ( 3n + 2) ( 6n + 5) ( 6n + 7)}{4( 4n + 1) {(4n + 3)}^2 (4n + 5)} 
H_n^2 
\end{equation}
Using the values
$H_0( \threequarters ) = 1$ and  
$H_1( \threehalves ) = \frac{7}{2}$, 
and calculating as in the derivation of 
(\ref{3over26}), we obtain 
\begin{equation}
\label{3over44}
H_n( \threequarters ) 
\newcommand{\seh}{{\textstyle\frac72}}
= (\seh)^n  \prod_{j=1}^{n-1} \prod_{i=1}^{n-j} \frac{9( 3i + 1) ( 3i + 2) ( 6i + 5) ( 6i + 7)}{4( 4i + 1) {(4i + 3)}^2 (4i + 5)} 
\end{equation}
The expression (\ref{3over44}) can be written in terms of factorials as given in (\ref{3over45}).
\begin{corollary}
\label{3over4corollary}
Suppose $a_k(x)$ is defined as in 
(\ref{akintro31})
and $H_n(x)= \det [ a_{i+j}(x) ]_{0 \leq i,j \leq n}$. Then
\begin{equation}
\label{3over45}
H_n( \threequarters ) 
= \prod_{i=1}^n \frac{(2i-1)! (6i+1)!}{2 (4i-1)! (4i+1)!}
\end{equation}
\end{corollary}

\section{A third identity and the differential equation for $H_n$}
\label{lastpiece}

We now continue with the proof of Theorem \ref{thm1}.
Using the expression for $ d_x^2 H_n$ in terms of $ H_n, K_n$ and $N_n$ from 
(\ref{SecondDerivativeAs013}),
and the expression for $ d_x H_n$ in terms of $ H_n$ and $K_n$ from 
(\ref{det01Eq31})
{\small
\begin{eqnarray}
\nonumber
&&(x-3) (2x-3)^2 (4x-3)^2  \Big( (x-1)(x-3)d_x^2 H_n  + (2(n+2)(x-3)+3)d_x H_n - 3n(n+1)H_n \Big) =
\\ 
\nonumber
&& 8(4n+3)(4n+5)(x-1)N_n  \\
\nonumber
&& - 4(4n+3)(-27n^2-93n-75 +(27n^2 + 81n+60)x +(8n+10)x^2 )K_n   \\
\label{rhsofde}
&&+ ( - 729n^4 - 3564 n^3 -6015n^2-4236n -1080+ (729n^4 + 2916n^3 + 3777n^2 + 1722n +180)x  + \\ 
\nonumber
&&\hspace*{2cm} (432n^3+ 1492 n^2 + 1676n +600 )x^2 ) H_n
\end{eqnarray}
}
To prove Theorem \ref{thm1}, we show that the expression on the right hand side of 
(\ref{rhsofde}) is identically zero. This is the consequence of our third identity.

Write the right hand side of (\ref{rhsofde}) in the form
\begin{equation}
\label{zero1}
- (p_n N_n + q_n K_n + r_n H_n )
\end{equation}
where
$p_n=p_n(x), q_n = q_n(x)$, and $r_n= r_n(x)$ 
are the negatives of the polynomials multiplying $ N_n$, $ K_n$, and $H_n$ on the right hand side of (\ref{rhsofde}),
respectively. 
That is,
{\small
\begin{eqnarray}
\nonumber
\hspace*{1cm} p_n(x)&=& -8(4n+3)(4n+5)(x-1) \\ \nonumber
\label{pnqnrn}
q_n (x)&=&  4(4n+3)(-27n^2-93n-75 +(27n^2 + 81n+60)x +(8n+10)x^2 )\\
r_n (x)&=& -( - 729n^4 - 3564 n^3 -6015n^2-4236n -1080+ (729n^4 + 2916n^3 + 3777n^2 \\ \nonumber
& & \hspace*{2cm} + 1722n +180)x  + (432n^3+ 1492 n^2 + 1676n +600 )x^2) 
\end{eqnarray}
}
Note that $p_n$ is linear, $ q_n$ and $r_n$ are quadratic polynomials in $x$.

The three matrices in (\ref{defHn}), (\ref{defKn}), (\ref{defNn}) 
that define  $H_n$, $K_n$, and $ N_n$ respectively, differ only in their last column. 
Therefore the expression
\begin{equation}
\label{zero2}
p_n N_n + q_n K_n + r_n H_n 
\end{equation}
is the determinant of the $ (n+1) \times (n+1)$  matrix whose 
\begin{enumerate}
\item[i)] first $n$ columns are the columns of the matrix $A_n$ in (\ref{defHn}),
\item[ii)] whose last column is the linear combination 
\end{enumerate}
\begin{equation}
\label{linear-comb}
p_nv_{n+2} + q_n v_{n+1} + r_n  v_n
\end{equation}
where for $ j\geq0$
$$v_j
        = \left[\begin{array}{c}
            a_j \\
            a_{j+1} \\
            \vdots \\
            a_{j+n}
                \end{array}\right] ~.
$$
To show that the determinant (\ref{zero2}) is zero,
it suffices to show that the last column (\ref{linear-comb}) of this matrix 
is a linear combination of the first $n$ columns $ v_0, v_1, \ldots, v_{n-1}$, i.e. there
are weights
$ w_{n,0}, w_{n,1}, \ldots, w_{n,n+2}$ with
\begin{eqnarray}
\nonumber
w_{n,n+2} & =&  p_n\\
\label{pqr}
w_{n,n+1}  &=& q_n\\
\nonumber
w_{n,n} &=&  r_n
\end{eqnarray}
such that
\begin{equation}
\label{weights}
\sum_{j=0}^{n+2} w_{n,j} a_{i+j} =0
\end{equation}
for $ i=0, 1, \ldots, n $. This is the third identity we need.

After some experimentation, bolstered by checks in Mathematica for up to $ n=30$, we found that the weights are likely given by
\begin{equation}
\label{weights2}
w_{n,j}(x) = \sum _{i=0}^{\min\left(-j+n+1,\left\lfloor
   \frac{n}{2}\right\rfloor +1\right)} 
A_{n,i,j}(x) ~~+ ~
\sum_{i=0}^{\min\left(-j+n+2,\left\lfloor \frac{n}{2}\right\rfloor +1\right)} 
B_{n,i,j} (x) 
\end{equation}
where
\begin{eqnarray*}
   A_{n,i,j}(x)  & =&  u_{n,i} ( \alpha_{n,i,j} + \beta_{n,i,j} x +  \gamma_{n,i,j} x^2)\\
   B_{n,i,j}(x)  &=& v_{n,i} ( \delta_{n,i,j} + \epsilon_{n,i,j} x +  \theta_{n,i,j} x^2)
\end{eqnarray*}
with
\dsp
{\small
\begin{eqnarray*}
    u_{n,i}  &=& 
\frac{3 \cdot 2^{i+1}(4n+3)(4n+5) 
 {n \choose 2i}
 {j+2(n-i+1)+1 \choose 2(n-i-j+1)+1}
(2n-i+2)!(4n-2i+5)!!(6n+7)!!}{(2i+1)(3j+1)(3j+2)(2i-j-2n-3)(2i-j-2n-2)(2n+1)!(4n+5)!!(6n-4i+7)!!}\\
 \alpha_{n,i,j}  &=& 
3(9j(j+1)+2)(12i^2+6(3j-4n-5)i+6(n+1)(2n+3)-j(18n+19)) \\
\beta_{n,i,j} &=&  
-3(3j+2)(24i^3+8(12j-9n-8)i^2+2(36n^2-96jn+64n+9j(5j-11)+23)i \\
&~&  \hspace*{3cm} -2(n+1)(2n+3)(6n+1)-j^2(90n+89)+3j(32n^2+66n+33)) \\
\gamma_{n,i,j}  &=& 
2(2i+2j-2n-3)(i+j-n-1)(36i^2+6(9j-12n-5)i+6n(6n+5)-3j(18n+13)+2) \\
    v_{n,i}  &=& 
\frac{3 \cdot 2^{i+3}(2n-2i-2j+5)(4n+3)(4n+5) 
 {n+1 \choose 2i}
 {j+2(n-i+2)+1 \choose 2(n-i-j+2)+1}
(2n-i+3)!(4n-2i+5)!!(6n+7)!!}{(3j+1)(3j+2)(2i-j-2n-5)(2i-j-2n-4)(2i-j-2n-3)(2n+2)!(4n+5)!!(6n-4i+9)!!} \\
 \delta_{n,i,j}  &=&  
3(9j(j+1)+2)(6i^2+3(3j-4n-7)i+3(n+2)(2n+3)-j(9n+14)) \\
 \epsilon_{n,i,j}  &=& 
-3(3j+2)(12i^3+(48j-36n-50)i^2+(45j^2-3(32n+49)j+4(n+1)(9n+16))i\\
&~&  \hspace*{3cm} -2(n+2)(2n+3)(3n+2)-j^2(45n+67)+3j(n(16n+49)+37)) \\
\theta_{n,i,j}  &=&  
2(2i+2j-2n-3)(i+j-n-2)(18n^2-9(4i+3j)n+33n-33j+3i(6i+9j-11)+13)
\end{eqnarray*}
}

\ssp

The proof of 
formula (\ref{weights}), and consequently the proof of Theorem \ref{thm1}, is 
therefore reduced to a somewhat involved system of binomial identities involving
summations over multiple indices. 
Our first reaction was that a pen-and-paper proof was infeasible. The explicit form of the conjectured weights gave
us a slim hope that automated binomial identity provers might work.
However, we did not see how to do this. So our entire approach to the evaluation of the determinants $H_n$
seemed to founder on these identities.
Our efforts seemed to indicate that 
these identities may well be  beyond the standard formulations and
computational resources required for automated binomial identity provers, and even if such a proof is possible,
the insight such a proof would offer would be minimal.

We were able to find an alternative approach, which avoids the explicit form of the weights, and which has the advantage
of generalization to other cases. For the proof, all 
we need to do is to show the {\em existence} of weights satisfying
(\ref{weights}) and (\ref{pqr}).
This becomes a manageable problem which relies on certain properties of the generating function 
\begin{equation}
\label{gf}
f = f(x,y) = \sum_{k \geq 0 } a_k(x) y^k
\end{equation}
of the $a_k$.

\begin{lemma}
\label{genfunc31}
Suppose $a_k(x)$ is as defined in (\ref{akintro31}) and $f$ is as in (\ref{gf}).
Then
\begin{equation}
\label{f3}
f(x,y) = \frac{4x+2t-6}{(x-3)(4x-3)y(2t-3) - (x-1)(27y-4)} 
\end{equation}
where 
$$t= \sum_{k\geq 0 }\frac{(3k)!}{(2k+1)! k!} y^k  = 1+ y +3y^2 +12y^3 + \cdots $$ satisfies $t^3 y = t-1$.
\end{lemma}
The proof can be found in Appendix I.

The proof of the existence of weights satisfying (\ref{weights}) and (\ref{pqr}) takes several steps.
Suppose $t$ is as in Lemma \ref{genfunc31} and define
\begin{equation}
\label{eta1}
\eta = 2t-3 = -1+ 2y + 6 y^2 + \cdots
\end{equation}
Then 
$$
f(x,y) = \frac{4x-3 + \eta}{(x-3)(4x-3)y\eta - (x-1)(27y-4)} 
$$
The first step is to note that we can find a nontrivial polynomial $Q_0 = Q_0(y)$ of degree
$n+1$ such that 
\begin{equation}
\label{2tm3}
\eta Q_0 = Q_1 + y^{2n+3} \Psi_0
\end{equation}
where $ Q_1= Q_1 (y)$ is a polynomial of degree $n+1$, and $ \Psi_0 = \Psi_0(y)$ is a power 
series in $y$. This is because there are $n+2$ coefficients in $Q_0$
and only $n+1$ linear equations these coefficients must satisfy.  
In the next step, put
\begin{equation}
\label{Q2}
Q_2 = (x-3) (4x-3) y Q_1 - (x-1) (27y-4) Q_0
\end{equation}
Then $Q_2= Q_2(x,y)$ is a polynomial in $x$ and $y$ of $ y$-degree $n+2$. It is important to
note that all three of the polynomials $Q_0, Q_1, Q_2$ so defined are nontrivial. We claim that the 
coefficients of $Q_2$ are the weights we want. In other words, the coefficients of the terms 
$y^{n+2}$ through $ y^{2n+2}$ in $ f Q_2$ vanish. Our aim now is to 
use the 
special form of the generating function $f$ in (\ref{f3}) to show this.

Writing 
(\ref{f3}) in the form
$$
f(x,y) \big( (x-3)(4x-3)y\eta - (x-1)(27y-4) \big) =
4x-3 +  \eta
$$
and multiplying through by $ Q_0$,  we get
$$
f(x,y) \Big( (x-3)(4x-3)y (Q_1 + y^{2n+3} \Psi_0) - (x-1)(27y-4)Q_0 \Big)
=
(4x-3) Q_0 + Q_1 + y^{2n+3} \Psi_0
$$
or 
$$
f(x,y) \Big( Q_2 + (x-3)(4x-3) y^{2n+4} \Psi_0  \Big)
=
(4x-3) Q_0 + Q_1 + y^{2n+3} \Psi_0
$$
and therefore
$$
f(x,y) Q_2 = 
(4x-3) Q_0 + Q_1 + y^{2n+3} \left( \Psi_0 + f(x,y) (x-3)(4x-3) y \Psi_0  \right)
$$
which means that
\begin{equation}
\label{fQ2}
f Q_2 =  (4x -3)Q_0 + Q_1 + y^{2n+3} \Psi_1
\end{equation}
where 
$$
(4x -3)Q_0 +Q_1
$$
is a polynomial in
$x$ and $y$ of $y$-degree $n+1$ and
$\Psi_1 = \Psi(y)$ is a power series in $y$. 

This last statement (\ref{fQ2}) is equivalent to 
the statement that 
\begin{equation}
\label{xxx}
\sum_{j=0}^{n+2} \mathcal{C}_{n+2-j}( Q_2) a_{i+j} ~=0
\end{equation}
for $ i = 0 ,1, \ldots , n$, where 
by $\mathcal{C}_k (\Psi)$ we denote the coefficient of the term $y^k$ in a power series $\Psi$.
Thus (\ref{weights}) holds with 
$$ w_{n,j} = \mathcal{C}_{n+2-j}( Q_2) $$ 
for $j = 0 ,1, \ldots , n+2$. 

In particular, taking determinants, we have that 
\begin{equation}
\label{ciden}
\mathcal{C}_{0}( Q_2)N_n + 
\mathcal{C}_{1}( Q_2) K_n +
\mathcal{C}_{2}( Q_2) H_n = 0
\end{equation}
This identity is not trivial. Otherwise, we would have a nontrivial linear relationship
among the $n+1$ columns $ v_0, v_1 , \ldots , v_n$ of $H_n$. But we have already shown that the special 
values $H_n(3), H_n ( \threehalves ) , H_n ( \threequarters )$ 
evaluated in (\ref{Cfactorials}), (\ref{3over27}), (\ref{3over45}) are nonzero, and therefore 
$H_n$ does not identically vanish. 

Now we need to verify the three special values in (\ref{pqr}). 
Since equation (\ref{xxx}) is valid when multiplied through by an arbitrary constant,
in terms of the coefficients $ \mathcal{C}_{k}( Q_2)$, we need to show that 
for some nonzero constant  $\alpha$,
\begin{eqnarray}
\nonumber
\mathcal{C}_{0}( Q_2) & =& \alpha  p_n\\
\label{cpqr}
\mathcal{C}_{1}( Q_2)  & = & \alpha q_n\\
\nonumber
\mathcal{C}_{2}( Q_2) &=& \alpha  r_n
\end{eqnarray}
where $p_n, q_n, r_n$ are as defined in (\ref{pnqnrn}).

First we rewrite (\ref{ciden}) in terms of 
the three coefficients
$\mathcal{C}_{0}( Q_0), \mathcal{C}_{1}(Q_0), \mathcal{C}_{2}( Q_0)$ which are pure 
constants, independent on $x$ and $y$.
First, comparing coefficients in 
(\ref{2tm3}) and (\ref{Q2}), we obtain the five identities
in (\ref{Csystem}).
\begin{eqnarray}
\nonumber
\mathcal{C}_{0}( Q_1)  & = & - \mathcal{C}_{0}( Q_0) \\ \nonumber
\mathcal{C}_{1}( Q_1)  & = & 2 \mathcal{C}_{0}( Q_0) - \mathcal{C}_{1}( Q_0) \\
\label{Csystem}
\mathcal{C}_{0}( Q_2) & =&  4(x-1) \mathcal{C}_{0}( Q_0) \\ \nonumber
\mathcal{C}_{1}( Q_2)  & = & (x-3)(4x-3)  \mathcal{C}_{0}( Q_1) -27(x-1) \mathcal{C}_{0}( Q_0)
+4(x-1) \mathcal{C}_{1}( Q_0) \\
\nonumber
\mathcal{C}_{2}( Q_2) &=&  (x-3) (4x-3)\mathcal{C}_{1}( Q_1) -27(x-1))\mathcal{C}_{1}( Q_0) +
4(x-1) \mathcal{C}_{2}( Q_0) 
\end{eqnarray}
Therefore
\begin{eqnarray}
\nonumber
\mathcal{C}_{0}( Q_2) & =&  4(x-1) \mathcal{C}_{0}( Q_0)\\
\label{ccpqr}
\mathcal{C}_{1}( Q_2)  & = & -2 (2x^2+6x-9) \mathcal{C}_{0}( Q_0) +4(x-1) \mathcal{C}_{1}( Q_0) \\
\nonumber
\mathcal{C}_{2}( Q_2) &=&  2(x-3)(4x-3) \mathcal{C}_{0}( Q_0)-
2(2x^2+6x-9)\mathcal{C}_{1}(Q_0)+ 4(x-1)\mathcal{C}_{2}( Q_0) 
\end{eqnarray}
and (\ref{ciden}) becomes
\begin{eqnarray}
\nonumber
& & 4(x-1) \mathcal{C}_{0}( Q_0) N_n + \\ 
\label{ciden2}
& & \left(-2 (2x^2+6x-9) \mathcal{C}_{0}( Q_0) +4(x-1) \mathcal{C}_{1}( Q_0) \right)K_n + \\  \nonumber
& & \left(2(x-3)(4x-3) \mathcal{C}_{0}( Q_0)- 2(2x^2+6x-9)\mathcal{C}_{1}(Q_0)+ 4(x-1)\mathcal{C}_{2}( Q_0)
\right) H_n = 0
\end{eqnarray}
We have explicit linear relationships between the three determinants $N_n, K_n, H_n$ at 
the points $ x = 3, \threehalves, \threequarters$.  
For $x=3$, we use the expressions in terms of $H_n(3)$ 
for $N_n(3)$ from  (\ref{x33}) and for $K_n(3)$ from (\ref{x33}) in (\ref{ciden2}) to obtain the linear
equation
\begin{eqnarray}
\nonumber
& &   8 (3 + 4 n) (5 + 4 n) \mathcal{C}_{2}( Q_0)  + \\ 
\label{ceq1}
& & \hspace*{5mm}6(5+4n)(1+2n+18n^2)  \mathcal{C}_{1}( Q_0) + \\  \nonumber
& & \hspace*{10mm}3(30+67n+338n^2+540n^3+243n^4) \mathcal{C}_{0}(Q_0)
= 0
\end{eqnarray}
For $x=\threehalves$, we use the expressions in terms of $H_n(\threehalves)$ 
for $N_n(\threehalves)$ from  (\ref{3over22}) and for $K_n(\threehalves)$ from (\ref{3over21}) in (\ref{ciden2}) to obtain 
\begin{eqnarray}
\nonumber
& & 8(3+4n)(5+4n) \mathcal{C}_{2}( Q_0)  + \\ 
\label{ceq2}
& & \hspace*{5mm}12(5+4n)(2+10n+9n^2) \mathcal{C}_{1}( Q_0) + \\  \nonumber
& & \hspace*{10mm}3(120+778n+1445n^2+1026n^3+243n^4) \mathcal{C}_{0}(Q_0)
= 0
\end{eqnarray}
Finally for $x=\threequarters$, we use 
(\ref{3over42}) and (\ref{3over41}) in (\ref{ciden2}) to obtain 
\begin{eqnarray}
\nonumber
& & 8(3+4n)(5+4n) \mathcal{C}_{2}( Q_0)  + \\ 
\label{ceq3}
& & \hspace*{5mm} 12(5+4n)(-4+n+9n^2) \mathcal{C}_{1}( Q_0) + \\  \nonumber
& &  \hspace*{10mm}3(-240-446n+95n^2+540n^3+243n^4) \mathcal{C}_{0}(Q_0)
= 0
\end{eqnarray}
Using (any two of) these equations 
we obtain the parametric solutions
\begin{eqnarray}
\nonumber
 \mathcal{C}_{2}( Q_0)  & = & 3(50n+343n^2+540n^3+243n^4) \alpha \\ 
\label{ceq4}
\mathcal{C}_{1}( Q_0) &= &  -12(1+n)(30+67n+36n^2) \alpha \\  \nonumber
\mathcal{C}_{0}(Q_0) & = & 8(15+32n+16n^2) \alpha
\end{eqnarray}
Going back to the system
(\ref{ccpqr}), we compute that 
\begin{eqnarray}
\nonumber
\hspace*{5mm} \mathcal{C}_{0}( Q_2) & = & 32(4n+3)(4n+5)(x-1) \alpha   \\
\label{cccpqr}
\mathcal{C}_{1}( Q_2)  & = & -16(4n+3)(-75-93n-27n^2+60x+81nx+27n^2x+10x^2+8nx^2) \alpha  \\
\nonumber
\mathcal{C}_{2}( Q_2) & = &  4(-1080-4236n-6015n^2-3564n^3-729n^4+180x+1722nx+ 3777n^2x+ \\ \nonumber
& & \hspace*{8mm} 2916n^3x+ 729n^4x+600x^2+1676nx^2+1492n^2x^2+432n^3x^2) \alpha
\end{eqnarray}
Taking $ \alpha = -\quarter$, 
these are exactly the weights $ w_{n,n+2} = p_n$, $w_{n,n+1}= q_n$, $w_{n,n}= r_n$ as claimed in (\ref{pqr}).
This finishes the proof of Theorem \ref{thm1}.
\end{proof}

\subsection{The differential equation solution for the $(3,1)$-case}

In the previous section we derived the differential equation 
(\ref{de1}) 
satisfied by 
$$
  H_n (x) = 
    \det \left[\begin{array}{cccccc}
          a_0 & a_1 & a_2 & \ldots & a_n \\
          a_1 & a_2 & a_3 & \ldots & a_{n+1} \\
       \vdots &     &     &        & \vdots \\
          a_n & a_{n+1}    & \ldots    &        & a_{2n}
       \end{array}\right]
$$
where
$$
       a_k = a_k^{(3,1)}(x)= \sum_{m=0}^k {3k+1-m \choose k-m} x^m.
$$
We obtain a power series solution to (\ref{de1})
of the form
$$
y (x) = \sum_{i=0}^{\infty}b_i (x-3)^i
$$
where the $b_i$ satisfy
$$
b_{i+1} = \frac{ (3 n +i +3) (n-i)}{(k+1)(2k+3)} b_i
$$
with $b_0 =1$.
This allows us to prove that 
\begin{equation}
\label{HandC}
  H_n (x) = C_n \sum_{i=0}^n\frac{n! (3n+i+2)!2^i(x-3)^i}
                                                   {(3n+2)! (n-i)!(2i+1)!}
\end{equation}
where $C_n$ is the constant of integration that 
depends only on $n$ but not on $x$.  
We have already calculated  
$ C_n = H_n (3)$
in (\ref{Cfactorials}) of Corollary \ref{3corollary}.
We record the final determinant evaluation result of the $(3,1)$-case as a Theorem:
\begin{theorem}
\label{final3np1}
Suppose $a_k(x)$ is defined as in 
(\ref{akintro31})
and $H_n(x)= \det [ a_{i+j}(x) ]_{0 \leq i,j \leq n}$. Then
\begin{equation}
\label{det3np1}
H_n (x)  =
      (-1)^n \prod_{i=1}^n \frac{ (6i-3)!(3i+2)!(2i-1)!} {(4i-1)!(4i+1)!(3i-2)!}
   \sum_{i=0}^n\frac{n!(3n+i+2)! 2^i(x-3)^i}
                    {(3n+2)! (n-i)! (2i+1)!}.
\end{equation}
\end{theorem}
In particular, evaluating (\ref{det3np1}) at $ x=0$ and at $ x=1$ we obtain
\begin{corollary}
\label{det3np1corollary}
{\small
\begin{eqnarray}
\det \left[  3 (i +j) +1 \choose i+j \right]_{0 \leq i,j \leq n} & =& 
(-1)^n \prod_{i=1}^n \frac{ (6i-3)!(3i+2)!(2i-1)!} {(4i-1)!(4i+1)!(3i-2)!} 
\sum_{i=0}^n\frac{n!(3n+i+2)! (-6)^i} {(3n+2)! (n-i)! (2i+1)!} \\
\det \left[  3 (i +j) +2  \choose i+j \right]_{0 \leq i,j \leq n} & =& 
(-1)^n \prod_{i=1}^n \frac{ (6i-3)!(3i+2)!(2i-1)!} {(4i-1)!(4i+1)!(3i-2)!} 
\sum_{i=0}^n\frac{n!(3n+i+2)!  (-4)^i} {(3n+2)! (n-i)! (2i+1)!} 
\label{second120}
\end{eqnarray}
}
\end{corollary} 

Note that the second evaluation (\ref{second120}) is an alternate form of the product formula (\ref{3np2}).

\subsection{A recursion for $H_n^{(3,1)} (x)$}
\label{3termsection}

We now remark on a property of the polynomials
$H_n(x)$ for the $(3,1)$-case, which does not enter into 
the evaluation of the determinant, but which nevertheless is of interest. 
\begin{theorem}
\label{3termrecursion}
With $H_0 = 1$, $H_1 (x)= -2x+5$,
the $H_n(x)$ satisfy a three-term polynomial recursion 
\begin{equation}
\label{3termHn}
p_n (x) H_n (x) + q_n(x) H_{n-1}(x) + H_{n-2}(x) =0
\end{equation}
for $ n \geq 2$ 
where 
\begin{eqnarray*}
p_n(x) &=& \frac{4 (4 n-3)!^2 (4 n-1)!^2}{9 (3 n-2)^2 (3 n-1)^2 (2 n-1)!^2 (6 n-5)!^2} (x-1)^2 \\
q_n(x) &=& \frac{4 (n-1) (4 n-5)! (4 n-1)!}{3 (3 n-2) (3 n-1) (4 n+1) (2 n-1)!  (6 n-5)!}
\Big( 8 (4 n-3) (4 n+1) x^3\\
&& \hspace*{1cm} -36 (4 n-3) (4 n+1) x^2 + 6 (126 n^2-63 n-23) x-3 (108 n^2-54 n-19) \Big)
\end{eqnarray*}
\end{theorem}
\begin{proof}
The theorem follows from the explicit form of 
the $H_n(x)$ given in Theorem \ref{final3np1}. We omit the details.
\end{proof}
The following property of the $H_n$ for the $(3,1)$-case is a consequence of Theorem \ref{3termrecursion}:
\begin{corollary}
The roots of 
$H_n(x)$ are real and interlace those of $H_{n+1}(x)$ for $ n>0$. 
\end{corollary}
\begin{proof}
The polynomials $H_n(x)$ are clearly not orthogonal, so the proof of the interlacing phenomenon does not follow
from classical results, nor directly from the results of \cite{Liu06}. 
Nevertheless an elementary induction argument works, which
we also omit.
\end{proof}

It is also curious that the recursion (\ref{3termHn}) seems to define a 
sequence of rational functions with nontrivial denominators. We are unaware of a general theory of 
such recurrences which actually define Sturm sequences of polynomials.

Another consequence of 
the recursion in 
Theorem \ref{3termrecursion} is a completely different derivation of the 
product formula (\ref{3np2}) for the $ (3,2)$-case. Letting $x=1$ in 
(\ref{3termHn}), and noting $p_n(1)$ vanishes,  we obtain the recursion
\begin{eqnarray*}
H_n^{(3,2)} & =&  \frac{-1}{q_{n+1}(1) } H_{n-1}^{(3,2)} \\
& = & \frac{3 (3 n+1) (3 n+2) (2 n+1)! (6 n+1)!}{4 n (4 n+1) (4 n-1)! (4 n+3)!} H_{n-1}^{(3,2)} 
\end{eqnarray*}
from which the product evaluation 
$$
 H_n^{(3,2)} = \prod_{i=1}^n \frac{3 (3 i+1) (3 i+2) (2 i+1)! (6 i+1)!}{4 i (4 i+1) (4 i-1)! (4 i+3)!}
$$
immediately follows.

\section{Transformation rules}
\label{lineartransformation}

Before we proceed with the evaluation of the determinant $ H_n ^{(2,1)} (x)$, we would like to 
collect the information we already have about various trace calculations to automate this process to some extent.
In the course of the calculations involved in the $(3,1)$-case, we have had to evaluate traces such as
$$
                 \mbox{Tr}(A_n^{-1}[p(x)a_{i+j+1}]_{0 \leq i,j\leq n})  ,~~
                 \mbox{Tr}(A_n^{-1}[q(x)(i+j)a_{i+j+1}]_{0 \leq i,j\leq n})  ,~~
                 \mbox{Tr}(A_n^{-1}[ r(x) c_{i+j}]_{0 \leq i,j \leq n})     
$$
for polynomial coefficients $p(x), q(x), r(x)$. These traces are reduced to calculating traces of the form
$$
                 \mbox{Tr}(A_n^{-1}[a_{i+j+1}])  ,~~
                 \mbox{Tr}(A_n^{-1}[(i+j)a_{i+j+1}])  ,~~
                 \mbox{Tr}(A_n^{-1}[ c_{i+j}])     
$$
etc., and then use linearity of the trace. Furthermore, the calculations do not 
use the explicit form of the polynomials
$a_k$. We can go one step further: when we apply the operator $\mbox{Tr}(A_n^{-1}*)$ to various terms of an identity
linearly, we can then multiply through the resulting transformed identity by $H_n$ to obtain 
a relationship between various determinants related to $H_n$, as we did
in a number of cases in the computation of the $(3,1)$-case. 
In Table \ref{lineartr} we summarize the effect of this transformation 
on various terms that may appear in an identity.

{\dsp
\begin{table}[htb]
\begin{center}
\begin{tabular}{|l|l|l|}
\hline 
$d_x a_n \rightarrow d_x H_n $& $~~$ \\
\hline \hline
& $ n a_{n-1} \rightarrow 0 $& $c_{n-1} \rightarrow 0$\\
\hline
$a_n \rightarrow (n+1)H_n $ & $n a_n \rightarrow n(n+1)H_n $& $c_n \rightarrow (2n+1)a_0 H_n $\\
\hline
$ a_{n+1} \rightarrow  K_n $ & $  n a_{n+1} \rightarrow 2nK_n $& $c_{n+1} \rightarrow 2 a_0 K_n + 2 n
a_1 H_n $\\
\hline
$a_{n+2} \rightarrow M_n + N_n~~~~ ( n >0 ) $& $ n a_{n+2} \rightarrow 2(n-1)M_n+ 2nN_n ~~~~ ( n >0 ) $
& \\
\hline 
\end{tabular}
\end{center}
\vspace*{5mm}
\caption{Some linear transformation rules.}
\label{lineartr}
\end{table}
}

In Table \ref{lineartr},
$ a_n = a_n(x), n = 0,1, \ldots $ is a sequence of polynomials, 
$c_n $ are the convolution polynomials defined by 
\begin{equation}
\label{convo}
c_n = \sum_{i=0}^n a_i a_{n-i}
\end{equation}
and $H_n$, $K_n$, $M_n$, $N_n$  are the 
following $(n+1) \times (n+1)$ determinants defined by the $a_k$'s:\\
$$
  H_n = \det \left[\begin{array}{cccccc}
          a_0 & a_1 & a_2 & \ldots & a_n \\
          a_1 & a_2 & a_3 & \ldots & a_{n+1} \\
       \vdots &     &     &        & \vdots \\
          a_n & a_{n+1}    & \ldots    &        & a_{2n}
       \end{array}\right]
~~
   K_n = \det \left[\begin{array}{cccccc}
       a_0     & a_1     & \ldots & a_{n-1}  & a_{n+1} \\
       a_1     & a_2     & \ldots & a_{n}    & a_{n+2} \\
       \vdots  &         & \ddots                        \\
       a_{n-1} & a_n     & \ldots & a_{2n-2} & a_{2n} \\
       a_{n} & a_{n+1} & \ldots & a_{2n-1}   & a_{2n+1}
              \end{array}\right]\\
$$

\ \\
$$
   M_n = \det \left[\begin{array}{ccccccc}
       a_0     & a_1     & \ldots & a_{n-2} &a_{n+1}  & a_{n} \\
       a_1     & a_2     & \ldots & a_{n-1} &a_{n+2}    & a_{n+1} \\
       \vdots  &         & \ddots                        \\
       a_{n-1} & a_n     & \ldots & a_{2n-3}& a_{2n} & a_{2n-1} \\
       a_{n} & a_{n+1} & \ldots & a_{2n-2 } & a_{2n+1}   & a_{2n}
              \end{array}\right]
 ~~
   N_n = \det \left[\begin{array}{ccccccc}
       a_0     & a_1     & \ldots & a_{n-2} & a_{n-1}  & a_{n+2} \\
       a_1     & a_2     & \ldots & a_{n-1} & a_{n}    & a_{n+3} \\
       \vdots  &         & \ddots  &                        \\
       a_{n-1} & a_n     & \ldots & a_{2n-3} & a_{2n-2} & a_{2n+1} \\
       a_{n} & a_{n+1} & \ldots & a_{2n-2} & a_{2n-1}   & a_{2n+2}
              \end{array}\right] \\
$$

From these rules we obtain a few known applications that are not connected with the 
determinants we are evaluating. We give a few simple examples.

\begin{example}
{\em Recall that a polynomial sequence
$ \{a_k\}_{k \geq 0}$ is said to be an {\em Appell set }  if $ d_xa_n = n a_{n-1}$. 
If $B_n$, $ {\mathcal H}_n$, $T_n$, $P_n$ are the Bernoulli, Hermite, Chebyshev, and Legendre polynomials,
respectively, then for any complex $ a \neq 0$, 
$$
\{ B_n (x) \}, ~~
\{ (2 a )^{-n} {\mathcal H}_n (ax) \}, ~~
\{ (x^2- a^2 )^{\frac{n}{2}} T_n (\frac{x}{\sqrt{x^2 -a^2}}) \}, ~~
\{ (x^2- a^2 )^{\frac{n}{2}} P_n (\frac{x}{\sqrt{x^2 -a^2}}) \}, 
$$
are all well-known examples of Appell sets \cite{Lupas88,Wimp00}.
Starting with the defining identity
$$
d_xa_n = n a_{n-1}
$$
and applying the transformations
$d_x a_n \rightarrow d_x H_n $ and  $ n a_{n-1} \rightarrow 0 $ from  Table \ref{lineartr},
we see that
$d_x H_n =0 $, and therefore the Hankel determinants of Appell polynomials are independent of $x$.}
\end{example}

\begin{example}
{\em For a sequence of polynomials $ a_n$ satisfying
\begin{equation}
\label{integrable}
x d_x a_n = n a_n +cn a_{n-1}
\end{equation}
for some constant $c$, 
we make the replacements
$d_x a_n \rightarrow d_x H_n $,  $ n a_{n} \rightarrow n(n+1)H_n $ and  $ n a_{n-1} \rightarrow 0 $ 
from Table \ref{lineartr}
to find that the corresponding Hankel determinant $H_n$ satisfies
$$
x d_x H_n = n(n+1) H_n 
$$
This is integrable with
$$
H_n(x) = c_n \cdot x^{n(n+1)}
$$
For $c=-1$ in (\ref{integrable}) we get the 
Laguerre polynomials 
$L_n$
with $L_0 = 1, L_1 =1 -x$. 
After simple row operations,
$c_n = H_n(1)$  has a known product evaluation as 
$$
c_n = \det \left[ \frac{1}{(i+j)!} \right]_{0 \leq i,j \leq n} =
 2^{-n^2} (-1)^{\frac{n (n + 1)}{2}}  \prod_{i=1}^n \frac{1}{1^2 \cdot 3^2 \cdots (2i-1)^2}
$$
For $c=1$ in (\ref{integrable}) we get the
so-called {\em derangement polynomials} 
$$
D_n(x) = \sum_{i=0}^n (-1)^i \frac{n!}{i!} x^{n-i}
$$
for which the constant of integration is 
$$
c_n = \prod_{i=1}^n i!^2
$$
as found in \cite{radoux98}.
}
\end{example}
\begin{example}
{\em  Legendre polynomials satisfy
$$
(x^2-1) d_x P_n = n x P_n - n P_{n-1}
$$
with $P_0 =1, P_1 =x $ and therefore
$$
(x^2-1) d_x H_n =  x n(n+1) H_n 
$$
and
$$
H_n = c_n \cdot (x^2-1)^{\frac{n(n+1)}{2}} ~.
$$
$H_n$ can be evaluated at $x=0$, giving
$$
c_n= 2^{-n^2} ~.
$$
}
\end{example}

\begin{example}
{\em 
In general, for any family of polynomials $ \{ a_k \}_{k \geq 0}$ such that 
$$
d_x a_n, ~~ a_n , ~~  n a_n, ~~ n a_{n-1}, ~~ c_n, ~~ c_{n-1}
$$
satisfy a linear identity  
with fixed polynomial coefficients
(where $c_n$ is the convolution defined in (\ref{convo})), the corresponding Hankel determinant
is the solution to a 
simple separable first-order linear equation.}
\end{example}

Now we return to the determinants $H_n^{(\beta, \alpha)}$ and consider the $(2,1)$-case. \\

\section{The $(2,1)$-case}

Now we let
\begin{eqnarray}\nonumber
a_k  & =&  a_k^{(2,1)} (x)  = \sum_{m=0}^{k}{2k+1-m \choose k-m}x^{m}\\
\label{2n1}
c_k &=& \sum_{i=0}^k a_i a_{k-i} \\ \nonumber
H_n & = & H_n^{(2,1)} (x) = \det  [a_{i+j}(x)]_{0\le i,j\le n}
\end{eqnarray}

\begin{theorem}
\label{defor2np1}
Let the polynomials $a_k = a_k^{(2,1)} (x)$ 
and $H_n= H_n^{(2,1)} (x)$ are as defined in 
(\ref{2n1}). Then
$H_n$ satisfies the differential equation
\begin{equation}
\label{de2np1}
x(x-2)d_x^2y + (2x-1)d_xy - n(n+1)y = 0
\end{equation}
\end{theorem}

The proof of Theorem \ref{defor2np1} mimics the steps of 
the $(3,1)$-case.  
The analogues of the first two identities in Lemma \ref{lemma1} and Lemma \ref{lemma2} of the $(3,1)$-case 
are now as follows:
\begin{lemma}
\label{thm12np1}
Suppose $ a_k$ and $ c_k$ are as defined in (\ref{2n1}). Then
\begin{equation}
\label{thm12np1proof}
2x(x-2)d_x a_n - (n+1)a_{n+1} + (4n+2(x+1))a_n - (x-1)c_n + 4(x-1 )c_{n-1} = 0
\end{equation}
\end{lemma}

\begin{lemma}
\label{thm22np1}
Suppose $ a_k$ and $c_k$ are as defined in (\ref{2n1}). Then
\begin{eqnarray}\nonumber
\label{thm22np1proof}
& & (nx+2(x+1))a_{n+2} - (2x(x+2)n + 2 (2x^2+3x +4)) a_{n+1}\\
& & \hspace*{2cm} + 4x^2(2n+3) a_n +(x-1)(x-2) c_{n+1} -4(x-1)(x-2) c_n = 0
\end{eqnarray}
\end{lemma}
The proofs are given in Appendix II.

Applying the transformations in Table \ref{lineartr} to the identities in these two Lemmas, 
we immediately obtain the pair of identities
\begin{eqnarray*}
&& 2x(x-2) d_x H_n - (2n+1)K_n + 4n(n+1)H_n +2(x+1)(n+1)H_n - (x-1)(2n+1)H_n = 0\\
&& x (2(n-1)M_n + 2n N_n) +2(x+1)(M_n+N_n) - 2x(x+2)(2 n K_n) \\
& & \hspace*{1cm} + 2 (2x^2+3x +4) K_n + 8x^2n(n+1)H_n +12x^2 (n+1)H_n \\
& &\hspace*{2cm}  +(x-1)(x-2) (2K_n+2n(x+3)H_n) -4(x-1)(x-2) (2n+1)H_n = 0
\end{eqnarray*}
These can be rewritten as 
\begin{eqnarray}\nonumber
\label{firstder2np1}
&& 2x(x-2) d_x H_n + (x+ 3+ 8n +4n^2)H_n-(2n+1)K_n =0 , \\ 
& &(2n x^3+ (8n^2+12n+8)x^2+(12n+10)x-4n-8)H_n \\ \nonumber
& & - ((2+4n)x^2+(12+8n)x +4)K_n + 2(n x+1)M_n +2((n+1)x +1)N_n =0
\end{eqnarray}
First we use identity 
(\ref{DeltaOneAsTrace31}) to rewrite $K_n$ as a trace in (\ref{firstder2np1}). 
Differentiating, we follow through
calculations of traces similar to the derivation of $d_x^2 H_n$ of the $ (3,1)$-case to obtain
\begin{eqnarray}
\nonumber
&  & 2(x-2)^2x^2(1+n x) d_x^2 H_n \\ \nonumber
& & \hspace*{5mm} + 2(x-2)x(3+10n +4n^2+3x+9n x+12n^2x+4n^3x+4nx^2+2n^2x^2) d_x H_n \\ \nonumber
& &  \hspace*{1cm} + (10+49n+78n^2+44n^3+8n^4+(4+34n+86n^2+98n^3+48n^4 + 8n^5)x\\
\label{d2for2np1}
& & \hspace*{3cm} +(1+8n+14n^2+8n^3)x^2+(2n-2n^3)x^3) H_n \\ \nonumber
& & \hspace*{4cm} -(2n+1)(1+2n+2n x+2n^2x) N_n  =0
\end{eqnarray}
Combining  
(\ref{d2for2np1}) and the expression (\ref{firstder2np1}) for $d_x H_n$, we obtain
\begin{eqnarray}
\label{3rn2np1}
& &  2 (x-2)x(nx+1) \Big( x(x-2)d_x^2y + (2x-1)d_xy - n(n+1)y \Big) = \\  \nonumber
& &  (2n+1)(1+2n+2nx+2n^2x) \left( N_n -  (4+2n+x )K_n +(2+5n+2n^2+3x+2nx)H_n \right)
\end{eqnarray}

\subsection{Product form evaluations at special $x$ for the $(2,1)$-case}

At this point, we can evaluate $H_n(x)$ at special points $x$ easily. At $x=0$, 
we get from 
(\ref{firstder2np1}) and  (\ref{d2for2np1}) 
\begin{eqnarray}
\nonumber
K_n &=& \frac{3+8n+4n^2}{2n+1} H_n\\ 
\label{2np1x0}
N_n &=& \frac{10+49n+78n^2+44n^3+8n^4}{(2n+1)^2}H_n \\ \nonumber
M_n &=& -(3n+2n^2 )H_n \nonumber
\end{eqnarray}
Using 
(\ref{detident})
and that  $H_0=H_1 = 1$ for $ x=0$, we obtain the recursion
for $\frac{H_n}{H_{n-1}}$ as
\begin{equation}
\label{2np1atx0}
\frac{H_{n+1}}{H_n}=  \frac{H_n}{H_{n-1}}  
\end{equation}
and 
therefore
\begin{equation}
\label{2np1atx0for}
H_n (0) =  1
\end{equation}
as is well known.  At $x=2$
\begin{eqnarray}
\nonumber
K_n &=& \frac{5+8n+4n^2}{2n+1}H_n\\ 
\label{2np1x2}
N_n &=& \frac{22+33n+20n^2+4n^3}{2n+1}H_n\\ \nonumber
M_n &=& -\frac{n(7+8n+4n^2)}{(2n+1)}H_n \nonumber
\end{eqnarray}
and $H_0=1$, $H_1 = -3$. From
(\ref{detident})
$$
H_{n+1} H_{n-1} = \frac{ (2n-1)(2n+3)}{(2n+1)^2}  H_n^2
$$
and therefore
\begin{equation}
\label{2np1constant}
H_n(2) = (-1)^n (2n+1)
\end{equation}

After the evaluations at special points, we now
return to the  proof of the differential equation for the $2n+1$ case. By (\ref{3rn2np1}), it suffices to 
prove the identity
\begin{equation}
\label{thirdiden2np1}
N_n -  (4+2n+x )K_n +(2+5n+2n^2+3x+2nx)H_n = 0
\end{equation}
This is the determinantal form of the third identity for the $(2,1)$-case.
The left hand side of (\ref{thirdiden2np1}) is a determinant as in 
(\ref{zero2}) and (\ref{linear-comb}) of the $(3,1)$-case.
The problem is again to show the existence of weights 
$ w_{n,0}, w_{n,1}, \ldots, w_{n,n+2}$ satisfying
the third identity (\ref{pqr}) and (\ref{weights}) where in this case
the requirements are
\begin{eqnarray}
\nonumber
w_{n,n+2} &=& p_n(x)= 1 \\ 
\label{pnqnrn2np1}
w_{n,n+1} &=& q_n (x)=   -(4+2n+x) \\ \nonumber
w_{n,n}&=& r_n (x)= 2+5n+2n^2+3x+2nx
\end{eqnarray}
By experimentation, we found that these weights are explicitly given by
\begin{equation}
\label{weights2np1}
w_{n,j}= (-1)^{n+j} \left\{ 2 {n+j+1 \choose n-j+1} + {n+j+1 \choose n-j+2} +
\left( 2 {n+j+1 \choose n-j} + {n+j+1 \choose n-j+1} \right) x \right\}
\end{equation}
for $ i=0, 1, \ldots, n-1$, and $w_{n,n}, w_{n,n+1}, w_{n,n+2}$ as in 
(\ref{pnqnrn2np1}).
The resulting identities are simpler than the $(3,1)$-case, and in fact a combination of 
paper-and-pencil and automatic identity provers can be used to prove that these are correct.
However we can do better, and we use this case as 
another exercise for the application of the method of the existence of the 
weights that we used in the $(3,1)$-case.

We need a lemma similar to Lemma \ref{genfunc31}. 
\begin{lemma}
\label{genfunc21}
Suppose $a_k(x)$ is as defined in (\ref{aks2np1}) and $f$ is as in (\ref{gf}).
Then
\begin{equation}
\label{f2np1}
f(x,y) = \frac{t}{(x-2) y t +1 - 2 x y}
\end{equation}
where
$$t= \sum_{k\geq 0 }\frac{(2k)!}{(k+1)! k!} y^k  = 1+ y +2y^2 +5y^3 + \cdots $$ satisfies $t^2 y = t-1$.
\end{lemma}
The proof of Lemma \ref{genfunc21} can be found in Appendix I.

Let $t$ be the power series in $y$ given in 
Lemma \ref{genfunc21}.
Then there is a nontrivial polynomial $Q_0 = Q_0(y)$ of degree
$n+1$ such that 
\begin{equation}
\label{2np12tm3}
t Q_0 = Q_1 + y^{2n+3} \Psi_0
\end{equation}
where $ Q_1= Q_1 (y)$ is a polynomial of degree $n+1$, and $ \Psi_0 = \Psi_0(y)$ is a power 
series in $y$; i.e. the coefficients of $y^k$ in $ t Q_0$ vanish for $ n+2 \leq k \leq 2n+2$. 
Such a nontrivial $Q_0$ exists 
because there are $n+2$ coefficients and only $n+1$ linear equations to satisfy. 
In the next step, put
\begin{equation}
\label{2np1Q2}
Q_2 =  (x-2) y Q_1+(1 - 2 x y)  Q_0
\end{equation}
Then $Q_2= Q_2(x,y)$ is a polynomial in $x$ and $y$ of $ y$-degree $n+2$.
All three polynomials $Q_0, Q_1, Q_2$ are nontrivial. We claim that the 
coefficients of $Q_2$ are the weights we want. In other words, the coefficients of the terms 
$y^{n+2}$ through $ y^{2n+2}$ in $ f Q_2$ vanish. 
Writing 
(\ref{f2np1}) in the form
$$
f(x,y) \big( (x-2) y t +1 - 2 x y  \big) = t
$$
and multiplying through by $ Q_0$,  we get
$$
f(x,y)  \Big( (x-2) y (Q_1 + y^{2n+3} \Psi_0) +(1 - 2 x y)  Q_0 \Big) = Q_1 + y^{2n+3} \Psi_0
$$
or 
$$
f(x,y) Q_2 =  Q_1 + y^{2n+3} \Psi_0 -f(x,y) (x-2)y y^{2n+3} \Psi_0 
$$
and therefore
$$
f(x,y) Q_2 = 
Q_1 + y^{2n+3} \left( \Psi_0 - f(x,y)  (x-2) y  \Psi_0  \right)
$$
which means that
\begin{equation}
\label{2np1fQ2}
f Q_2 =  Q_1 + y^{2n+3} \Psi_1
\end{equation}
where
$\Psi_1 = \Psi(y)$ is a power series in $y$. 
This last statement (\ref{2np1fQ2}) is equivalent to 
\begin{equation}
\label{2np1xxx}
\sum_{j=0}^{n+2} \mathcal{C}_{n+2-j}( Q_2) a_{i+j} ~=0
\end{equation}
for $ i = 0 ,1, \ldots , n$, where 
by $\mathcal{C}_k (\Psi)$ we denote the coefficient of the term $y^k$ in a power series $\Psi$.
Thus (\ref{weights}) holds with 
$$ w_{n,j} = \mathcal{C}_{n+2-j}( Q_2) $$ 
for $j = 0 ,1, \ldots , n+2$. 

Therefore
\begin{equation}
\label{2np1ciden}
\mathcal{C}_{0}( Q_2)N_n + 
\mathcal{C}_{1}( Q_2) K_n +
\mathcal{C}_{2}( Q_2) H_n = 0
\end{equation}
This identity is not trivial, for otherwise we would have a nontrivial linear relationship
among the $n+1$ columns $ v_0, v_1 , \ldots , v_n$ of $H_n$. But we have already shown by the evaluations in 
(\ref{2np1atx0for}), (\ref{2np1constant}) that $H_n$ does not identically vanish. 

Now we need verify the three special values in (\ref{pnqnrn2np1}), i.e.
for some nonzero constant  $\alpha$,
\begin{eqnarray}
\nonumber
\mathcal{C}_{0}( Q_2) & =& \alpha  p_n\\
\label{2np1cpqr}
\mathcal{C}_{1}( Q_2)  & = & \alpha q_n\\
\nonumber
\mathcal{C}_{2}( Q_2) &=& \alpha  r_n
\end{eqnarray}
where $p_n, q_n, r_n$ are as defined in (\ref{pnqnrn2np1}).

Rewrite (\ref{2np1ciden}) in terms of 
$\mathcal{C}_{0}( Q_0), \mathcal{C}_{1}(Q_0), \mathcal{C}_{2}( Q_0)$ which are pure 
constants, independent on $x$ and $y$.
Comparing coefficients in 
(\ref{2np12tm3}) and (\ref{2np1Q2}), we obtain 
\begin{eqnarray}
\nonumber
\mathcal{C}_{0}( Q_1)  & = & \mathcal{C}_{0}( Q_0)  \\ \nonumber
\mathcal{C}_{1}( Q_1)  & = &  \mathcal{C}_{1}( Q_0) + \mathcal{C}_{0}( Q_0) \\ 
\label{2np1Csystem}
\mathcal{C}_{0}( Q_2)  & = &  \mathcal{C}_{0}( Q_0) \\ \nonumber
\mathcal{C}_{1}( Q_2)  & = & (x-2) \mathcal{C}_{0}( Q_1) - 2 x \mathcal{C}_{0}( Q_0) +\mathcal{C}_{1}( Q_0)  \\
\nonumber
\mathcal{C}_{2}( Q_2) &=&  (x-2) \mathcal{C}_{1}( Q_1) - 2x \mathcal{C}_{1}( Q_0) + \mathcal{C}_{2}( Q_0)
\end{eqnarray}
Therefore
\begin{eqnarray}
\nonumber
\mathcal{C}_{0}( Q_2) & =&  \mathcal{C}_{0}( Q_0)\\
\label{2np1ccpqr}
\mathcal{C}_{1}( Q_2)  & = & \mathcal{C}_{1}( Q_0) -(x+2) \mathcal{C}_{0}( Q_0) \\
\nonumber
\mathcal{C}_{2}( Q_2) &=&  \mathcal{C}_{2}( Q_0)- (x+2)\mathcal{C}_{1}(Q_0) +(x-2) \mathcal{C}_{0}( Q_0) 
\end{eqnarray}
and (\ref{2np1ciden}) becomes
\begin{eqnarray}
\nonumber
& &  \mathcal{C}_{0}( Q_0) N_n + \\ 
\label{2np1ciden2}
& &\hspace*{15mm}  \left( \mathcal{C}_{1}( Q_0) -(x+2) \mathcal{C}_{0}( Q_0)  \right) K_n + \\  \nonumber
& & \hspace*{30mm}  \left( \mathcal{C}_{2}( Q_0)- (x+2)\mathcal{C}_{1}(Q_0) +(x-2) \mathcal{C}_{0}( Q_0)  \right) H_n = 0
\end{eqnarray}
Using the expressions for $K_n$ and $N_n$ in terms of $H_n$ at $x=0$ and $x=2$ from 
(\ref{2np1x0}) and (\ref{2np1x2})
\begin{eqnarray*}
& & (10+49n+78n^2+44n^3+8n^4)\mathcal{C}_{0}( Q_0) \\
& & \hspace*{15mm} +(3+8n+4n^2)(2n+1)(\mathcal{C}_{1}( Q_0) - 2 \mathcal{C}_{0}( Q_0)) \\
& & \hspace*{30mm} + (2n+1)^2 ( \mathcal{C}_{2}( Q_0)- 2 \mathcal{C}_{1}(Q_0) -2 \mathcal{C}_{0}( Q_0)) =0 \\
& & (22+33n+20n^2+4n^3)\mathcal{C}_{0}( Q_0) \\
& & \hspace*{15mm} + (5+8n+4n^2)(\mathcal{C}_{1}( Q_0) -4 \mathcal{C}_{0}( Q_0)) \\
& & \hspace*{30mm} + (2n+1)( \mathcal{C}_{2}( Q_0)- 4 \mathcal{C}_{1}(Q_0)) = 0
\end{eqnarray*}

The parametric solutions are
\begin{eqnarray}
\nonumber
 \mathcal{C}_{2}( Q_0)  & = & n(2n+1) \alpha \\ 
\label{2np1ceq4}
\mathcal{C}_{1}( Q_0) &= &  -2(n+1) \alpha \\  \nonumber
\mathcal{C}_{0}(Q_0) & = &  \alpha
\end{eqnarray}
Going back to the system (\ref{2np1ccpqr})
\begin{eqnarray}
\nonumber
\hspace*{5mm} \mathcal{C}_{0}( Q_2) & = & \alpha   \\
\label{2np1cccpqr}
\mathcal{C}_{1}( Q_2)  & = & -(4+2n+x) \alpha  \\
\nonumber
\mathcal{C}_{2}( Q_2) & = & (2+5n+2n^2+3x+2nx) \alpha 
\end{eqnarray}
With $ \alpha = 1 $, 
these are exactly the weights $ w_{n,n+2} = p_n$, $w_{n,n+1}= q_n$, $w_{n,n}= r_n$ as claimed in
(\ref{pnqnrn2np1}).

\subsection{The differential equation solution for the $(2,1)$-case}
\label{21desection}

We have the differential equation (\ref{de2np1})
satisfied by the Hankel determinants $H_n$ for the $(2,1)$-case.
We obtain a power series solution to (\ref{de2np1})
as
$$
y(x) = C_n \sum_{i=0}^{n}\frac{(n+i)!2^i (x-2)^i}{(n-i)! (2i+1)!}
$$
where he constant of integration $C_n = H_n (2)$ is given by 
(\ref{2np1constant}). Therefore
\begin{theorem}
\label{final2np1}
Suppose $a_k(x)$ is defined as in 
(\ref{aks2np1})
and $H_n(x)= H_n^{(2,1)}(x)= \det [ a_{i+j}(x) ]_{0 \leq i,j \leq n}$. Then
\begin{equation}
\label{det2np1}
H_n( x ) =    (-1)^n (2n+1) \sum_{i=0}^{n}\frac{(n+i)!2^i (x-2)^i}{(n-i)! (2i+1)!}
\end{equation}
\end{theorem}

Evaluating (\ref{det2np1}) at $x=0$ and at $ x=1$ we obtain
\begin{corollary}
\label{2np1cor}
\begin{eqnarray}
\label{2np2x0}
\det \left[  2 (i +j) +1  \choose i+j \right]_{0 \leq i,j \leq n}  & =&
 (-1)^n (2n+1) \sum_{i=0}^{n}\frac{(n+i)!(-4)^i}{(n-i)! (2i+1)!}\\
\label{2np2}
\det \left[  2 (i +j) +2  \choose i+j \right]_{0 \leq i,j \leq n}  & = &
 (-1)^n (2n+1) \sum_{i=0}^{n}\frac{(n+i)!(-2)^i}{(n-i)! (2i+1)!}
\end{eqnarray}
\end{corollary}
Interesting point about Corollary \ref{2np1cor} is that 
the right hand side of 
(\ref{2np2x0}) is a complicated way of writing 1, while
the right hand side of (\ref{2np2}) evaluates to the simple expression 
\begin{equation}
\label{alternate}
(-1)^\frac{n(n+1)}{2} ~.
\end{equation}
Alternate derivations of the evaluation (\ref{alternate}) for the  determinant in 
(\ref{2np2}) were communicated to us by Ira Gessel and Christian Krattenthaler.
A generalization of the determinant in (\ref{2np2}) can be found in 
Corollary \ref{2np2jacobi2} below.

It is also interesting that the determinants
$ H_n^{(2,1)}(x)$ are orthogonal polynomials:

\begin{corollary}
\label{2np2jacobi1}
The $(n+1) \times (n+1)$ Hankel determinants 
$ H_n (x)=  H_n^{(2,1)}(x)$ form an orthogonal family of polynomials.
\end{corollary}
\begin{proof}
We have 
$$ H_0(x) = 1, ~~ H_1 (x) = 1- 2x 
$$
and for $ n \geq 2$,
\begin{equation}
\label{3term2np1}
H_n(x) = 2(1-x) H_{n-1}(x) - H_{n-2}(x)
\end{equation}
which can be verified directly from the explicit formula in 
(\ref{det2np1}) or using the differential equation 
(\ref{de2np1}).
Orthogonality now follows from Favard's theorem.
\end{proof} 

\noindent
In this case the generating function
of the $H_n(x)$ is
$$
\frac{ 1-y} { 1 - 2 (1-x)y + y^2} = \sum_{n\geq 0} H_n(x) y^n ~.
$$
Comparing with the 
generating function 
$$
\frac{ 1} { 1 - 2 xy + y^2} 
= \sum_{n\geq 0} U_n(x) y^n 
$$
of the Chebyshev polynomials of the second kind $U_n(x)$, 
we have
\begin{equation}
\label{Chebyshev}
H_n(x)= U_n(1-x) - U_{n-1}(1-x) ~.
\end{equation}
In fact we can prove more. 
From the identity
\begin{equation}
\sum_{m =0}^k {2k+2 \choose  k-m} x^m =  \sum_{m = 0}^k {2  k + 1 - m \choose k - m}(x + 1)^m = a_k^{(2,1)}(x+1) 
\end{equation}
and the differential equation
(\ref{de2np1}) for $H_n^{(2,1)}(x)$, 
we obtain both the evaluation of the Hankel determinant of the polynomials
\begin{equation}
\label{2np2ak}
a_k(x) = \sum_{m =0}^k {2k+2 \choose  k-m} x^m ~,
\end{equation}
and an alternate expression for $H_n^{(2,1)}(x)$ itself.
\begin{corollary}
\label{2np2jacobi2}
Suppose the polynomials $ a_k(x)$ are as defined in 
(\ref{2np2ak}) and $H_n (x) = \det [ a_{i+j} (x) ]_{0 \leq i,j \leq n}$. Then
$H_n(x)$ satisfies Jacobi's differential equation
\begin{equation}
\label{de2np2ak}
(x^2-1)d_x^2y + (2x+1)d_xy - n(n+1)y = 0
\end{equation}
and therefore
in terms of the Jacobi polynomials
\begin{eqnarray*}
H_n(x)  & = &  (-4)^n \frac{(n!)^2}{(2 n)!} P_n^{( \half , - \half)} (x), \\
H_n^{(2,1)}(x) & = & (-4)^n \frac{(n!)^2}{(2 n)!} P_n^{( \half , - \half)} (x-1)~.
\end{eqnarray*}
\end{corollary}

\section{Remarks}
\label{conclusions}

When we started to look at Hankel determinants of $A_n = [ a_{i+j}]_{0 \leq i,j \leq n}$ 
with entries $a_k = {3k+1 \choose k}$, it looked
like an almost product could account for the combination of small and
large prime factors that were showing up in the data.  At the time, we
also thought that whatever was moving the $(3,1)$-case off the pure
product formulas of the $(3,0)$ and $(3,2)$-cases would also explain the
$(3,3)$-case, and the $(3,4)$-case, and so on. It might also get at
the $(4,1)$-case, and all the rest of the $(\beta , \alpha)$-cases.

This section is a brief discussion of what did and did not work out.
We started with the observation that if you define $ a_k(x)$ by 
(\ref{akintro31})
then 
$$
a_k(1) = {3k+2 \choose k} ,
$$
so that 
$$
H_n^{(3,2)} = H_n^{(3,1)} (1)~
$$ has a product form evaluation.  What we wanted was the value of 
$H_n^{(3,1)}(x)$ not at $x=1$ but at $x=0$.  
Experimenting with the expansion of 
$H_n^{(3,1)} (x)$ at various points, we found that 
at $ x=1$ and at $x=3$
the Taylor coefficients of $H_n^{(3,1)} (x)$ 
had small prime factors.
We called $x=1$ and $x=3$ {\em round points}.

We also found several Hankel determinants that had
similar behavior but we had trouble finding a consistent
explanation for this phenomenon.  
In others, such as 
$H_n^{(3,0)}$,  the
coefficients were somewhat round but not round enough to allow us to
find a simple power series expansion.

When we tried to determine the common elements of these different
determinants  differential equations came to mind because
hypergeometric series such as~(\ref{mainformula}) clearly have a
differential  equation.  In retrospect, the explanation of why $x=1$ and $x=3$ are a
good places 
to expand $H_n^{(3,1)}$ is partially explained by examining the form
of its differential equation.  For instance (\ref{de1}) can be written in the form
\begin{equation}
\label{operator}
(x-3)^2d_x^2y +2(n+2)(x-3)d_xy- 3n(n+1)y = - 2(x-3)d_x^2y - 3 d_xy.
\end{equation}
This form is significant because the operator on the left hand side of (\ref{operator})
takes $(x-3)^k$ to $(x-3)^k$ times a polynomial in $n$ and $k$ and the
operator on the right 
takes $(x-3)^k$ to $(x-3)^{k-1}$ times a polynomial in $n$ and $k$.
This means that the differential equation defines a two-term recursion
on successive coefficients of the expansion of $H_n^{(3,1)}(x)$ around
$(x-3)$.  So the power series is easily shown to be hypergeometric.

Differential equations can also explain the ``near'' round behavior seen 
in some cases, such as the $H_n^{(3,0)}$.
However,  the existence of a differential equation would  not provide any input on why the
constant of integration such as $H^{(3,0)}_n(3)$ was round. Consequently it was a nice
bonus when, during the proof of the differential equation, we realized
that we could prove product formulas for these kinds of determinants
(e.g., (\ref{Cfactorials}) in Corollary \ref{3corollary},
(\ref{3over27}) in Corollary \ref{3over2corollary}, 
(\ref{3over45}) in Corollary \ref{3over4corollary};
(\ref{special1a}), 
(\ref{special1b}), and
(\ref{special2})
in Section \ref{additional}).

Once differential equations entered the picture, it was easy to guess the
differential equations for many of these Hankel determinants and obtain
strong experimental evidence for their
correctness, but we did not find a general
pattern, even for the families
$H_n^{(\beta, \alpha)}$. 
It is doubtful that such a pattern exists.
For example,
the (unproven) differential equation for the $(3,2)$-case is
fourth order and quite mysteriously complex (see
Figure~\ref{3k+2diffeq}).  This differential equation has been tested
for $H_n^{(3,2)}(x)$ where $n=1, 2, \ldots, 75$.  We did not find a
differential equation for $(3,3)$ or for any $(3,\alpha)$ with
$\alpha\ge 3$.

\begin{figure}
$$
  \begin{array}{l@{}l}
 2&(3x-1)^2(x-3)^3x(4(n+2)(2n+1)x^2+(8n^2+20n+11)x-1)\;d_x^4y \\
  &{} + (x-3)^2(3x-1)\begin{array}[t]{l@{}l}
           (12&(n+2)(2n+1)(8n+27)x^4 
                  - 3(128n^3+568n^2+724n+161)x^3 \\
              &{} - (576n^3+3016n^2+4756n+2269)x^2 
                  + (72n^2+252n+319)x - 15)\;d_x^3y
                       \end{array} \\
  &{}+3(x-3)\begin{array}[t]{l@{}l}
         (12&(8n^4+118n^3+427n^2+533n+174)x^5 \\
            &{} -(736n^4+6704+19628n^2+21289n+5814)x^4 \\
            &{}+(800n^4+2944n^3+564n^2-7580n-7078)x^3 \\
            &{}+2(816n^4+6744n^3+19358n^2+23069n+9809)x^2 \\
            &{}-6(108n^3+540n^2+972n+679)x+15(9n+20))d_x^2y
            \end{array} \\
  &{}-3\begin{array}[t]{l@{}l}
         (12&(16n^5+62n^4-7n^3-293n^2-378n-120)x^5 \\
            &{}-(960n^5+3104n^4-3080n^3-20993n^2-23419n-6450)x^4 \\
            &{}+4(144n^5-4n^4-2300n^3-4657n^2-2040n+838)x^3 \\
            &{}+(1728n^5+6624n^4-5832n^3-556990n^2-80626n-36004)x^2 \\
            &{}+12(318n^4+1995n^3+4670n^2+4928n+2058)x -(783n^2+3105n+3102))d_xy
        \end{array} \\
  &{}-3n(n+1)\begin{array}[t]{l@{}l}
              (12&(12n^4+68n^3+137n^2+113n+30)x^4 \\
                 &{} -(864n^4+4488n^3+8158n^2+5887n+1222)x^3 \\
                 &{}+(720n^4+2304n^3+280n^2-4547n-3388)x^2 \\
                 &{}+3(576n^4+3816n^3+9182n^2+9533n+3666)x 
                    -3(120n^2+507n+538))y = 0
             \end{array}  \\
  \end{array}
$$
\caption{Differential equation for $H_n^{(3,2)}(x) $. \label{3k+2diffeq}}
\end{figure}

Now the actual proof of the differential equation for a determinant 
like $ H_n^{(3,1)} (x)$ is an elusive creature. The straightforward evaluation of
the derivative yields a sum of $n+1$ determinants, each one of which is
badly behaved at $ x=3$. That was a disturbing development, and we still had
to differentiate each of these determinants themselves to get at the second
order differential equation.

In order to approach the problem systematically, we wrote the
derivative of a Hankel determinant  
$H_n = \det [ a_{i+j} (x) ]_{0 \leq i,j \leq n}$ as 
$$
     d_x H_n = \mbox{Tr}(A_n^{-1}d_xA_n)H_n
$$
where
$$
     d_xA_n =
     d_xA_n (x) =
    [d_x a_{i+j}(x)]_{0\le i,j\le n} .
$$
The obvious thing to try at this point is to look for expansions of the form
$$
    d_x A_n = R A_n + A_n S
$$
so that we could prove that 
$$
    d_x \det(A_n) = \mbox{Tr}(R + S) \det(A_n).
$$
This search also proved fruitless, and in retrospect Lemma~\ref{firstder}
suggests that it will be very difficult to find $R$ 
and $S$ because
$$
   (x-3)(2x-3)(4x-3)\mbox{Tr}(R + S) =  
           \left( 8nx^2-6(5n+1)x-3(9n^2+13n+8) \right)  +2(4n+3) \frac{K_n}{H_n}
$$
and  experiments suggest that $K_n$ and $H_n$ are often relatively prime.

So you have to feed some carefully selected facts into $d_x A_n$ to
get a useful expansion of
$ d_x H_n$.  Initially the first order 
differential equation
$$
x (x-1) d_x a_n(x)  -(n (x-3) -2) a_n(x)  -(3n+2) a_n(0) = 0
$$
looks like the fact needed.  However this does not work because we
never figured out what to do with the evaluation of the term
$$
   \mbox{Tr}(A_n^{-1}[(3i+3j+2)a_{i+j}(0)]_{0 \leq i,j \leq n}).
$$

So at this point we started searching for matrices that behaved well
after being evaluated by the operator  
\begin{equation}\label{troperatorForRemarks}
    \mbox{Tr}(A_n^{-1}*)
\end{equation}
By ``behaves well'' we mean that the calculations return a single determinant,
or at most a linear combination of a just a few determinants, thus avoiding the expansion of $H_n$ as a sum of 
a large number of determinants.  We quickly learned about the first
two columns in Table~\ref{lineartr} in Section \ref{lineartransformation}.  But
this was not  
sufficient because $d_xA_n$ cannot be expressed as a sum of the
matrices defined in these two columns. 

At this point we recalled some results in the
literature~\cite{Iwasaki02,Kajiwara01} that related derivatives to
convolutions.  We realized that what this work told us was that the
operator in (\ref{troperatorForRemarks}) works out very nicely on
$[c_{i+j+1}]$, $[c_{i+j}]$, $[c_{i+j-1}]$, i.e. the trace came out to
be a single determinant as already shown in Table \ref{lineartr}.

This gave us enough tools to find our first identity, given as
(\ref{FirstId31}) in Lemma \ref{lemma1}.  Using this identity we have
enough to express the first and second derivatives of $H_n$ as linear
combinations of small numbers of determinants.  In particular in the
examples of this paper, the first derivative of $H_n$ is expressed as
a linear combination of $H_n$ and $K_n$.  The second derivative is
expressed as a linear combination of $H_n$, $K_n$, $M_n$ and $N_n$.
This process can easily be continued to higher derivatives though the
computations get messier at each derivative.

If the first identity (Lemma~\ref{lemma1}) is weakened to include terms like $a_{n+2}$
and 
$c_{n+2}$ the above process still works.  The only difference is that the
derivatives have more distinct determinant summands.  The point of the above
discussion is that we cannot handle any $ (\beta , \alpha)$-case in which
there does not exist a first identity of this kind.  So far there are a dozen
cases, some not of the $ (\beta , \alpha)$ type, which have a first identity
that we can handle, and we mention $ (3,0)$, $(3,2)$, and $(2,2)$ as examples
of these.

The role of the second (Lemma~\ref{lemma2}) and third identities
(~\ref{weights}) is to prove linear relationships
between the determinants that are generated by the above process of
differentiating $H_n$. 

Of these three identities, the proof is most sensitive to the form of the
second identity.  
In the $(3,2)$-case, the second identity involves $a_{n+3}$ terms.  
While this version of the second identity can still be used to prove linear
relationships between determinants, there are not enough linear relationships,
and a proof of the fourth order differential equation  for the 
$(3,2)$-case does not seem possible with the tools of this paper. 

Another difficulty that needs to be overcome is the problem of finding and 
proving the third identity.  Our original approach was to guess the form of
the third identity and thus reduce the problem to binomial identity proving
techniques.  This works quite well for the $(2,1)$-case where the third
identity has an explicit form:

\begin{equation}
\sum_{i=0}^{n+2}(-1)^i\left(2{n+1+i\choose n+1-i}+{n+1+i\choose n+2-i} +
                   \left(2{n+1+i\choose n-i}+{n+1+i\choose n+1-i}\right)x\right)a_{i+m}(x)=0
\end{equation}
for $0\le m\le n$.

Though we were eventually able to guess the third identity in the
$(3,1)$-case, its form was hardly enlightening.  In addition, the identity
was so complex that even the job of applying automated tools to prove 
it would be a major undertaking. So the
whole process had arrived at an impasse.

The way forward was a sidestep via an existence theorem
(Section~\ref{lastpiece}).  The
third identity (once proved) said that the binomial identities we
wished to prove actually {\em existed}, and were unique, and you only
had to know a few of their components explicitly to show that the
right hand side of (\ref{rhsofde}) vanished. The third identity 
for the $ (3,1)$-case 
is (\ref{weights}).  This identity is proved via the generating function
of the $ a_n(x)$ and from the explicit linear relationships governing
the determinants $ H_n(x)$, $K_n(x)$, $N_n(x)$ at special values of
$x$.

Of special mention is the form of the generating function given in
Lemma \ref{genfunc31}. It is important to know that this is not the
form that emerges in the straightforward derivation of the generating
function, as was done in Lemma \ref{genfunc}.  We were surprised that
we could not use this form directly to prove the third identity.  It
took a long time to discover the form in Lemma \ref{genfunc31}. This
latter form worked, and this concluded the proof of the $(3,1)$-case. \\

To sum up, a lot of things have to fall right in place perfectly to
evaluate a $H_n^{(\beta, \alpha )} (x)$.

\section{Additional results on Hankel determinants}
\label{additional}

The proof technique presented here is applicable to Hankel determinants of polynomials other than the 
$ a_k^{( \beta, \alpha)}$.
Along with the $(3,0)$-case, 
we give the necessary ingredients, i.e. the three identities required, 
for the proof of the differential equation satisfied for
a few of these but omit the proofs of the theorems and 
the construction of the explicit power series 
solutions.

\begin{example}
\label{example3m0}

{\em We start with the $(3,0)$-case

First identity for the $(3,0)$-case:
\begin{eqnarray*}
& & 3(x-3)x(4x-3)d_x a_n - (4(2x-3)n + 2(2x-5)) a_{n+1} \\
& & +  (27(2x-3) n + 3(4x^2-3x-9))a_n-(x-1)c_{n+1} + 27(x-1) c_n = 0
\end{eqnarray*}

Second identity for the $(3,0)$-case:
\begin{eqnarray*}
& & (4(2x -3)^2(5x-3)n + 2(2x-3)(5x-3)(6x-11))a_{n+2} \\
& & - (81 (8 x^3 - 24x^2 + 27x - 9)n  + 18 (37 x^3 - 123x^2 + 153x - 54)) a_{n+1} \\
& & + (729 x^3 n + 486x^3)a_n + 4(x-1)(2x-3)(5x-3) c_{n+2} \\
& & - 3 (40 x^4 - 30x^3 - 207x^2 + 270x - 81)c_{n+1} + 162 x^2(5 x^2 - 15x + 9)  c_n =0
\end{eqnarray*}

The generating function of the $ a_k$ for the third identity for the $(3,0)$-case:
$$
f(x,y)= \frac{-(2x-3)t-3x}{(x^2(9y-4)+10x-6)t + (x-3)(4x-3)}
$$
where $ t^3y = t-1$.

\begin{theorem}
\label{thme1}
The Hankel determinant for the 
$(3,0)$-case satisfies the differential equation
\begin{eqnarray*}
&&  (x-3)(2x-3)(5x-3)d_x^2y  - 2(10(n-1)x^2 - 9(3n-4)x-9(n+5))d_x y \\
&& \hspace*{2cm} + n(10(n-1)x - 3(n-7))y =0 ~.
\end{eqnarray*}
\end{theorem}

In the $(3,0)$-case, in addition to the already known product
form at $ x=0$ given in (\ref{3np0}), at $ x=3$ and 
surprisingly also at  $ x =\threehalves $ the determinant is given by a simple
product. We omit the proofs but record these evaluations below:

\begin{equation}
\label{special1a}
    H_n^{(3,0)}(3)=\frac{(3n)!(3n+2)!}{2 (n!^2)}
                      \prod_{i=1}^n\frac{3(6i-5)!(2i)!(2i-1)}
                                        {(4i+1)!(4i-1)!} ~,
\end{equation}
\begin{equation}
\label{special1b}
    H_n^{(3,0)} (\threehalves ) = \prod_{i=1}^{n}\frac{27(6i-5)!(3i-1)(3i-2)(2i-1)!} {2(4i-1)!(4i-3)(4i-4)!} ~.
\end{equation}
}
\end{example}

\begin{example}

{\em
Next, take 
\begin{equation}
\label{3k2m}
a_k(x) =   \sum_{m=0}^{k}{3k-2m \choose k-m}x^{m}
\end{equation}

The first identity for the polynomials in (\ref{3k2m}):
\begin{eqnarray*}
&& x(2x-9)(4x+9)d_x a_n - (36n + 30) a_{n+1}  + (243 n + 8x^2 + 18x + 81) a_n \\
&& - 12 c_{n+1}- (8x^2-36x-81) c_n + 27x(2x-9) c_{n-1} = 0
\end{eqnarray*}

The second identity for the polynomials in (\ref{3k2m}):
\begin{eqnarray*}
&& (36(2 x + 3) n + 66 (2 x+3)) a_{n+2} - ((32 x^3 +486x + 729) n + 12(4x^3 + 4 x^2 + 54x + 81)) a_{n+1} \\
& & + (216 x^3n  + 108 x^2(2  x + 3)) a_n + 12(2 x + 3) c_{n+2} + (16x^3-72 x^2 - 378 x-243) c_{n+1}\\
& & +  2 x (8 x^3-90x^2 + 243 x +729) c_n - 54  x^3 (2x-9) c_{n-1} =0
\end{eqnarray*}

The generating function of the $ a_k$ for the third identity for the polynomials in (\ref{3k2m}):
$$
f(x,y) =
\frac{3t + 2x}{(x(4 x-9) y - 6) t + 9 + 2x - 6x^2 y}
$$
where $ t^3y = t-1$.

\begin{theorem}
\label{thme2}
The Hankel determinant of the polynomials in (\ref{3k2m}) satisfies the 
differential equation 
$$
(2x + 3) (2  x - 9) d_x^2 y + 4(2 (n + 2) x - 9 (n + 1)) d_x y - 12 n(n + 1) y = 0 ~.
$$
\end{theorem}
}
\end{example}

\begin{example}
\label{exampleaex}

{\em
What we refer to as the ``aex''-case (for exceptional) is the Hankel determinant where
\begin{equation}
\label{aex}
a_k(x) =   \frac{1}{k+1} \sum_{m=0}^{k}{3k-m \choose k-m} (m+1)(m+2) x^{m}
\end{equation}

The first identity for the aex-case:
\begin{eqnarray*}
&&x (x - 3)(7x-3) d_x a_n - (4(x - 1)n + 14 (x - 1)) a_{n+1} \\
& & + (27 (x - 1) n+3 (x - 1) (11x + 3)) a_n - 6 (x - 1)^2 c_n + 6 x^3 c_{n-1} = 0
\end{eqnarray*}

The second identity for the aex-case:
\begin{eqnarray*}
&& (4(x-1)^2(3x - 1) n +  18 (x - 1)^2(3x - 1)) a_{n+2} + ((-113x^3 + 189 x^2 -135 x + 27) n \\
& &  - 4(30 x^4 - 19 x^3 - 21 x^2 + 39 x - 9)) a_{n+1}+ (216 x^3 n + 12 x^3(7 x^2 - 2 x + 15)) a_n\\
& & + 6 (x - 1)^3 (3 x - 1) c_{n+1}  - 6  x^2(x - 1)(10 x^2 -17x+9) c_n + 6 x^5(7x - 9) c_{n-1} = 0
\end{eqnarray*}

The generating function of the $ a_k$ for the third identity for the aex-case:
$$
f(x,y)= \frac{2 (x - 1)^2 \tau - 2 (x - 1)^2 - 2 (3 x - 1)}{(-2yx^2 (x - 3) \tau + 2 yx^2(x - 3) + 6 x^2y - 3 x + 1}
$$
where
$$
\tau = \sum_{k\geq 0 }\frac{(3k)!}{(2k)! (k+1)!} y^k  ~.
$$

\begin{theorem}
\label{thme3}
The Hankel determinant of the polynomials in (\ref{aex}) satisfies the 
differential equation 
$$
(3x-1)(x-1)(x-3) d_x^2 y - 2(3nx^2 - 8 nx - 3 (n + 4)) d_x y + 3 n(n + 1)(x - 1) y =0 ~.
$$
\end{theorem}

In the aex-case, 
the determinant is given by a simple
product for 
$ x =\threeoverseven $ as shown below:
\begin{equation}
\label{special2}
    H_n (\threeoverseven ) = \prod_{i=0}^n\frac{2(6i+7)!(2i+1)!}{7(4i+5)!(4i+3)!}
\end{equation}
}
\end{example}

\begin{example}
\label{example3k0m}

{\em
We can aslo evaluate the Hankel determinants as an almost product for 
\begin{equation}
\label{example3k0meq}
   a_k(x) = \sum_{m=0}^k {3k+1 \choose k-m}x^m.
\end{equation}

These polynomials are related to the polynomials $a_k^{(3,0)}(x)$ of Example \ref{example3m0} 
by the simple transformation
\begin{equation}
\sum_{m =0}^k {3k+1 \choose  k-m} x^m =  \sum_{m = 0}^k {3  k  - m \choose k - m}(x + 1)^m = a_k^{(3,0)}(x+1) 
\end{equation}
Therefore
\begin{theorem}
\label{thme4}
The Hankel determinant of the polynomials in (\ref{example3k0meq}) satisfies the
differential equation
{\small
$$
(x-2)(2x-1)(5x+2)d_x^2y - 2(10(n-1)x^2-(7n-16)x-26n-19)d_xy + n(10(n-1)x+7n+11)y = 0  ~.
$$
}
\end{theorem}
}
\end{example}

\begin{example}

{\em
Finally we have an alternate evaluation of $H_n^{(3,1)}(x)$ at $x=1$.
\begin{theorem}
\label{thme5}
Suppose $a_k(x)$ is defined as in 
(\ref{akintro31})
and $H_n(x)= \det [ a_{i+j}(x) ]_{0 \leq i,j \leq n}$. Then

\begin{equation}
\label{det3np1at1}
H_n(x)  =
 \prod_{i=1}^n\frac{(6i+4)! (2i+1)!}{2 (4 i+2)!(4i+3)!}
       \sum_{i=0}^n \frac{(-1)^i n!(4n+3)!!  (3n+i+2)! (x-1)^i}
                              {(3n+2)! i!(n-i)!(4n+2i+3)!!}.
\end{equation}
\end{theorem}

As in Corollary 
\ref{det3np1corollary},  taking $ x=0$ and $ x=1$ 
in (\ref{det3np1at1}) we obtain
the following alternate evaluations. 
\begin{corollary}
\label{det3np1at1corollary}
{\small
\begin{eqnarray*}
\det \left[  3 (i +j) +1 \choose i+j \right]_{0 \leq i,j \leq n} & =&
 \prod_{i=1}^n\frac{(6i+4)! (2i+1)!}{2 (4 i+2)!(4i+3)!}
       \sum_{i=0}^n \frac{ n!(4n+3)!!  (3n+i+2)! }
                              {(3n+2)! i!(n-i)!(4n+2i+3)!!}\\
\det \left[  3 (i +j) +2  \choose i+j \right]_{0 \leq i,j \leq n} & =&
 \prod_{i=1}^n\frac{(6i+4)! (2i+1)!}{2 (4 i+2)!(4i+3)!}
\end{eqnarray*}
}
\end{corollary}
The first evaluation in Corollary \ref{det3np1at1corollary} is 
the expression given in (\ref{mainformula2}) for $H_n^{(3,1)}$.
It is special as there are no cancellations on the right.
The second evaluation is identical in form to 
the known product in (\ref{3np2}) for $H_n^{(3,2)}$.
}
\end{example}

\ \\

There are other variants on the polynomials $ a_k^{(\beta, \alpha)}$ 
which experiments suggest satisfy differential equations.  
A particularly unusual example are the Hankel determinants for 
$$
a_k(x) = \sum_{m=0}^{k+1} { 3k + 4 -m \choose k+1-m}x^m
$$
These determinants a third order equation, but 
the coefficients are very large and not round, making it hard to guess what
they are.

Going back to the Hankel determinants $ H_n^{(\beta, \alpha)}$,
some of these are governed by a 
second order differential equation, such as
$(3,0)$, $(3,1)$, and $(2,1)$.
There are also non-$( \beta, \alpha)$-cases  governed by a
second order differential equation such as the aex-case in Example \ref{exampleaex}. 

In forthcoming work, we plan to develop the methods of this paper to investigate 
the Hankel determinants related to other families of polynomials. In the 
family $(2,r)$, for example, the three 
cases $ r= 0, 1, 2$  have product evaluations, but for $ r \geq 3$, the 
evaluations are no longer products but only almost products. We 
also encounter new phenomena in these cases, such as
higher order differential equations and 
case-splitting into 
residue classes.

We remark that experimentally we 
know that both $(2,3)$ and $(3,2)$-cases satisfy  fourth order differential equations.
For the $(2,3)$-case, 
the polynomial in front of the fourth derivative has  degree 11.

\ \\

\noindent
{\bf Acknowledgments:} We would like to thank Julius Borcea, Ira Gessel, Christian Krattenthaler,  Boris
Shapiro, and Doron Zeilberger. 

\newpage
\section{Appendix I: The generating function of the $a_k^{(\beta, \alpha)}(x)$}

In Lemma \ref{genfunc},
we give a closed form of the generating function $f(x,y) $ defined in (\ref{gf}).
\begin{lemma}
\label{genfunc}
Suppose $a_k(x)$ is as defined in (\ref{akintro}) and $f$ is as in (\ref{gf}).
Then
\begin{equation}
\label{fab}
f(x,y)  = \frac{t^{\alpha+1}}{(\beta + ( 1-\beta) t) (1-xyt^{\beta -1})}
\end{equation}
where 
\begin{equation}
\label{tcube}
 t^\beta y = t-1~.
\end{equation}
\end{lemma}
\begin{proof}
Changing the order of summation and rearranging
$$ f(x,y) = \sum_{m \geq 0 } x^m y^m \sum_{n \geq 0 } { \beta n + (\beta-1)m + \alpha \choose n} y^n $$
It is known \cite{PS1925} that
\begin{equation}
\label{polya}
\sum_{n \geq 0 } {   \alpha + \beta n \choose n} y^n  = \frac{t^{\alpha+1}}{ \beta + (1-\beta ) t }
\end{equation}
where  $t$ satisfies (\ref{tcube}).
For our generating function, $ \beta $ is the same but $ \alpha$ in (\ref{polya}) is replaced by 
$ ( \beta -1) m + \alpha $.
Using these parameters we obtain
\begin{equation}
\label{gfagain}
f(x,y)  = \sum_{m \geq 0 } x^m y^m \frac{t^{( \beta -1)m + \alpha +1 }}{ \beta+ (1 -\beta) t} ~=~
\frac{t^{\alpha+1}}{(\beta + ( 1-\beta) t) (1-xyt^{\beta -1})}
\end{equation}
where
$t$ satisfies (\ref{tcube}).
\end{proof}

Note that using the Lagrange inversion formula, $t$ can be expanded as 
\begin{equation}
\label{expandt}
t= \sum_{k\geq 0 }\frac{(\beta k)!}{((\beta -1)k+1)! k!} y^k = 1 + y + \beta y^2 + \frac{ \beta (3 \beta-1)}{2}
y^3 + \cdots 
\end{equation}

\noindent
{\bf Proof of Lemma \ref{genfunc31}:}
\begin{proof}
The expression (\ref{gfagain}) for $f$  with $ \beta = 3$ and $ \alpha =1$ 
can be rewritten in the form (\ref{f3}) since
$$
\frac{t^2}{(3-2t)(1-xyt^2)} -\frac{4x+2t-6 }{(x-3)(4x-3)y(2t-3) - (x-1)(27y-4)}
$$
is equal to
$$
\frac{(t-1-t^3 y)(2t x + 6x -9)}{(3-2t)(1-xyt^2 )(4 tx^2 y -6x^2 y-15txy+9xy+9ty+2x-2)} \\
$$
and the numerator has $t -1-t^3 y $ as factor, but  this is zero by 
(\ref{tcube}).  This gives the form of the generating function $f$ for the $ (3,1)$-case that is claimed in (\ref{f3}).
\end{proof}

\ \\
\noindent
{\bf Proof of Lemma \ref{genfunc21}:}
\begin{proof}
The proof is similar to the proof of
Lemma \ref{genfunc31}
once we observe that the expression for the generating function
\begin{equation}
f(x,y) = \frac{t^2}{(2-t) (1-x y t)}
\end{equation}
from Lemma \ref{genfunc} can be written as (\ref{f2np1}) since
$$
\frac{t^2}{(2-t) (1-x y t)} -
\frac{t}{(x-2)yt + 1-2xy}
=
\frac{2t ( 1-t+t^2  y)}{(t-2)(1-xyt)(1-2yt - 2xy - xyt)}
$$
and the right hand side vanishes since in this case $ \beta = 2$ and
\begin{equation}
\label{t2y}
1-t+ t^2 y  = 0 
\end{equation}
by Lemma \ref{genfunc}.
\end{proof}

\newpage
\section{Appendix II: Pairs of identities for the $(3,1)$ and $(2,1)$-cases. }

Proofs of the first two identities in the $(3,1)$-case are  as follows:

\noindent
{\bf Proof of Lemma \ref{lemma1}:}
\begin{proof}
We make use of the generating function $f=f(x,y)$ of the $a_k$'s in the form
given by (\ref{gfagain}) in Appendix I.
Passing to the generating functions, (\ref{FirstId31}) is equivalent to the functional identity
\begin{eqnarray}
\label{thm1proof}
&&   (x-3)(2x-3)(4x-3)d_xf -  4yd_y \frac{f-1}{y}-6 \frac{f-1}{y} +(8x^2-18x+36)f \\ \nonumber
& & \hspace*{2cm} +27y d_yf -4(2x^2-6x+3) f^2   +27(2x^2-6x+3)y f^2 =0
\end{eqnarray}
Using identity (\ref{tcube}), 
\begin{equation}
\label{deroft}
 d_yt = \frac{t^3}{1-3yt^2} ~.
\end{equation}
Substituting this expression in the computation of
$d_yf$, the left hand side of 
(\ref{thm1proof}) can be simplified as
\begin{eqnarray*}
&& \frac{2 (t^3 y -t +1)}{(1-3yt^2)}
\Big( 2 x t^2+12 x^2 y t^2-90 x y t^2+54 y t^2+2 t^2+8 x^2 t\\
& & \hspace*{2cm} +x t-36 x^2 y t+135 x y t-81 y t-15 t-8 x^2-3 x+9 \Big)
\end{eqnarray*}
which vanishes since $t^3y -t+1=0$ by (\ref{tcube}).
\end{proof}

\noindent
{\bf Proof of Lemma \ref{lemma2}:}
\begin{proof}
We again use the generating function 
$f=f(x,y)$ of the $a_k$'s in the form
given by (\ref{gfagain}).
The identity (\ref{SecondId31}) is equivalent to 
\begin{eqnarray} \nonumber
& & 8(x-1)yd_y\frac{f-1-(4+x)y}{y^2} + 20(x-1)\frac{f-1-(4+x)y}{y^2}  \\ \nonumber
& & -2(-81+135x-72x^2+16x^3) yd_y\frac{f-1}{y}  -2(-117+180x-92x^2+24x^3)\frac{f-1}{y}\\
\label{thm2proof}
& & +27(2x-3)^3 yd_yf +54(2x-3)(3-4x+2x^2)f +  8(x-1)(3-6x+2x^2))\frac{f^2-1}{y} \\ \nonumber
& & +2(81-297x+324x^2-114x^3+8x^4)f^2 -27x(2x-3)(9-12x+2x^2)yf^2 = 0
\end{eqnarray}
Using the expression in (\ref{deroft}) for $ d_yt $ in the calculation of $ d_yf$, the left hand side of 
(\ref{thm2proof}) can be simplified as
\begin{eqnarray*}
&  & \frac{2 (t^3 y -t +1)}{(3-2t)^2 y^2 (1-3yt^2)(1-xyt^2)^2}
\Big(
-72t^2y^2x^4+16t^2yx^4+56tyx^4-72yx^4-108t^2y^2x^3 \\
& & -16tx^3-48t^2yx^3 -156tyx^3+  252yx^3+16x^3-4t^2x^2 +1134t^2y^2x^2 +14tx^2\\
& & + 72t^2yx^2+252tyx^2-486yx^2-10x^2 -1701t^2y^2x+32tx+18t^2yx\\
& & -405tyx+567yx-24x+4t^2 +729t^2y^2-30t-54t^2y+243ty-243y+18 \Big)
\end{eqnarray*}
which vanishes since $t^3y -t+1=0$ by (\ref{tcube}).
\end{proof}

Proofs of the first two identities in the $(2,1)$-case are  as follows:

\noindent
{\bf Proof of Lemma \ref{thm12np1}:}
\begin{proof}
Passing to the generating functions in (\ref{thm12np1proof}), we need to prove the identity
$$
2x(x-2)d_xf - d_yf + 4yd_yf+2(x+1)f - (x-1)f^2 + 4(x-1)y f^2 =0.
$$
Using the expression 
\begin{equation}
\label{dyt2}
d_yt = \frac{t^2}{1-2yt} 
\end{equation}
in the calculation of $ d_yf$, the left hand side of 
the identity we want to prove
can be simplified as
$$
\frac{ 2 t^2 (t^2y-t+1)(2-t+2x+4ty-6txy)}{(2-t)^2 (1-2ty)(1-txy)^2 }
$$
which vanishes by (\ref{t2y}).
\end{proof}

The proof of the second identity for the $(2,1)$-case is as follows:

\ \\
\noindent
{\bf Proof of Lemma \ref{thm22np1}:}
\begin{proof}
Again passing to the generating functions in (\ref{thm22np1proof}), we need to prove the identity
\begin{eqnarray*}
& & xy d_y\frac{f-1-(x+3)y}{y^2} + 2(x+1)\frac{f-1-(x+3)y}{y^2} - 2x(x+2)y d_y\frac{f-1}{y}
-2 (2 x^2+3x+4) \frac{f-1}{y} \\
& & \hspace*{2cm} + ~8 x^2 y d_yf + 12 x^2 f + (x-1)(x-2)\frac{f^2-1}{y}-4(x-1)(x-2)f^2 = 0
\end{eqnarray*}
Using the expression for $f$ and the expression for $d_y t$ in (\ref{dyt2}), the left hand side of this 
expression can be simplified to
\begin{eqnarray*}
& & \frac{2( t^2y -t+1 ) }{{( 2 - t) }^2y^2 ( 1 - 2ty ) {( 1 - txy ) }^2}
    \Big(  t^2 -4 + 8ty - 4t^2y - 2t^3y + 8txy - 2t^2xy + t^3xy \\
& & \hspace*{2cm} - 2t^2x^2y   + 8t^3y^2 - 16t^2xy^2 - 4t^3xy^2 + 8t^2x^2y^2 + 4t^3x^2y^2 - 4t^3x^3y^3 \Big) 
\end{eqnarray*}
which again vanishes by (\ref{t2y}).
\end{proof}

\newpage
\section{Appendix III: On the degree of a class of Hankel Determinants}

\begin{theorem}\label{Main}
Let $p_0, p_1, \ldots,p_{n}$ and $q_0, q_1, \ldots,q_{n}$ be integer sequences.
Further let $\gamma$ be real and
$\alpha_0,\alpha_1,\ldots, \alpha_{n}$
and
$\beta_0,\beta_1,\ldots , \beta_{n}$
be sequences of real numbers.  Then the determinant
\begin{equation}
\label{detdegree}
  \det \left[
	\sum_{0\le m \le p_i+q_j}
 	   {\alpha_i+\beta_j+\gamma m \choose p_i+q_j-m}x^m
			  \right]_{0 \leq i,j \leq n}
\end{equation}
as a polynomial in $x$ has degree $ \leq \max \{ \max{p_i} + \max{q_j} - n, 0 \}$.
\end{theorem}

\begin{proof}
We note the convention that the empty sum is zero, and the binomial coefficient in 
(\ref{detdegree}) is interpreted via the gamma function.
The proof of the theorem is by induction.  At each stage of the
induction we make use of the following lemma.  

\begin{lemma}\label{EasyLemma}
Let $p_0,p_1,\ldots, p_n$ and $q_0,q_1,\ldots,q_n$
be two integer sequences and let $a_0,a_1,a_2,\ldots$ be an
infinite sequence of real numbers.  Then the determinant
$$
	\det \left[
	     \sum_{0\le m\le p_i+q_j}a_{p_i+q_j-m}x^m\right]_{0 \leq i,j \leq n}
$$
as a polynomial in $x$ has degree $ \leq \max \{ \max{p_i} + \max{q_j} - n, 0 \} $.
\end{lemma}
\begin{proof}
By rearranging the indices, we can assume that the sequences
$p_0, p_1, \ldots,p_{n}$
and
$q_0, q_1, \ldots,q_{n}$
are nondecreasing without changing the conclusion of the lemma.
Then
$$
	\det \left[
	     \sum_{0\le m\le p_i+q_j}a_{p_i+q_j-m}x^m\right]_{0 \leq i,j \leq n}
				=
		y^{-\sum p_i - \sum q_j}\cdot
	\det \left[
	     \sum_{0\le m\le p_i+q_j}a_m y^m\right]_{0 \leq i,j \leq n}
$$
where $xy=1$.  Therefore the analysis can focus on the degree of
$$
     \det \left[
	     \sum_{0\le m\le p_i+q_j}a_m y^m\right]_{0 \leq i,j \leq n}
$$
as a polynomial in $y$.  By elementary row operations, 
we see that this is equal to the following
determinant
$$
  \det \left[\left\{
	\begin{array}{ll}
	   {\displaystyle \sum_{0\le m\le p_0+q_j}a_m y^m}      &\mbox{if }i=0 \\
	   {\displaystyle \sum_{0 \leq m;  ~ p_{i-1}+q_j+1\le m\le p_i+q_j }
			a_m y^m }                &\mbox{otherwise}
	\end{array}
			   \right.\right]_{0 \leq i,j \leq n}.
$$
Whenever $ p_{n}+q_{n} \geq n$,
the highest power of
$y$ in this determinant is no more than
$$
	\max \{ \sum p_i + \sum q_j , 0 \} 
$$
and the lowest power of $y$ in this determinant is at least 
$$
	\max \{ n+ \sum_{i=0}^{n-1}p_i + \sum_{j=0}^{n-1} q_j  , 0 \}
$$
Now if we multiply this polynomial by $y^{-\sum p_i-\sum q_j}$ obtaining
$$
		y^{-\sum p_i - \sum q_j}
	\det \left[
	     \sum_{0\le m\le p_i+q_j}a_m y^m\right]_{0 \leq i,j \leq n}
$$
we get a polynomial in $y^{-1}$.  The 
degree of this polynomial is no more than 
$$
	\max \{ p_{n}+q_{n} - n, 0 \}.
$$
Converting back to $x$ and remembering that $p_i$ and $q_j$ are
nondecreasing gives the desired result.
\end{proof}

Now we are ready to start the proof of Theorem~\ref{Main}.  Let
\begin{eqnarray*}
   \alpha_i' & = & \alpha_i + \gamma p_i		\\
   \beta_j'  & = & \beta_j + \gamma q_j
\end{eqnarray*}
With this change of variable, 
$$
	{\alpha_i+\beta_j+\gamma m \choose p_i+q_j-m} = 
        {\alpha_i'+\beta_j'-\gamma (p_i+q_j-m) \choose p_i+q_j-m}.
$$
The following identity holds:
$$
  {\alpha_i+\beta_j+\gamma m \choose p_i+q_j-m}
= \sum_{r\ge 0} \sum_{s= 0}^{p_i+q_j-m - r}
	 {\alpha_i' \choose r}{\beta_j' \choose s}
		{-\gamma (p_i+q_j-m) \choose p_i-r+q_j-s-m}.
$$
This allows us to express the determinant as
$$
   \det \left[
      \sum_{r\ge 0} \sum_{s\ge 0}
         {\alpha_i' \choose r}{\beta_j' \choose s}
          \sum_{0\le m\le p_i-r+q_j-s}
		{-\gamma (p_i+q_j-m) \choose p_i-r+q_j-s-m} x^m
		           \right]_{0 \leq i,j \leq n}
$$
This determinant has an expansion as a sum of terms of the form
$$
         \left(\prod_i {\alpha_i' \choose r_i}\right)
	 \left(\prod_j {\beta_j' \choose s_j}\right).
  \det \left[
          \sum_{0\le m\le p_i-r_i+q_j-s_j}
	    {-\gamma p_i+ -\gamma q_j + \gamma m \choose p_i-r_i+q_j-s_j-m}
                     x^m
		           \right]_{0 \leq i,j \leq n}
$$
over collections of nonnegative integers 
$ r_0,r_1,\ldots , r_{n} $ and $ s_0,s_1,\ldots , s_{n} $.

If all the $r_i$ and $s_j$ in any such collection are zero, then the corresponding determinant is
as described in Lemma~\ref{EasyLemma} and therefore has the
appropriate degree as a polynomial in $x$.  If any of the $r_i$ or
$s_j$ are nonzero, then the determinant has the same form as that of
Theorem~\ref{Main} except that $p_i$ has been replaced with $p_i'= p_i-r_i$
and $q_j'= q_j-s_j$.  This is the induction process.

Note that the induction process terminates, since eventually 
$$ 
\max \{ \max{p_i'} + \max{q_j'} - n, 0 \} = 0
$$
and in fact the resulting determinants in (\ref{detdegree}) become zero.
\end{proof}

\begin{corollary}
The $(n+1) \times (n+1)$ Hankel determinant $H_n^{(\beta, \alpha )}$ defined by the sequence of polynomials
in (\ref{form}) is of degree at most $n$.
\end{corollary}
\begin{proof}
$H_n^{(\beta, \alpha )}$ is of the  type described in Theorem \ref{Main}  
where $p_i = q_i = i $, $ \alpha_i = \beta i + \alpha$,
$ \beta_i = \beta i $, and 
$\gamma = -1$. 
\end{proof}

\bibliographystyle{plain}

\end{document}